\documentclass[11pt]{amsart}
\usepackage{amscd,amsmath,amssymb,amsfonts,appendix}
\usepackage[cmtip, all]{xy}
\usepackage{enumerate}
\usepackage[usenames,dvipsnames]{xcolor}
\usepackage{url}
\usepackage{mathtools}

\usepackage[T2A,T1]{fontenc}
\newcommand{\Zhe}{\mbox{\usefont{T2A}{\rmdefault}{m}{n}\CYRZH}}
\newcommand{\Yu}{\mbox{\usefont{T2A}{\rmdefault}{m}{n}\cyryu}}

\usepackage[english]{babel}

\usepackage[OT2,OT1]{fontenc}
\newcommand\cyr
{\renewcommand\rmdefault{wncyr}
\renewcommand\sfdefault{wncyss}
\renewcommand\encodingdefault{OT2}
\normalfont
\selectfont
}
\DeclareTextFontCommand{\textcyr}{\cyr}

\usepackage[latin1]{inputenc}
\usepackage{times}
\theoremstyle{plain}
\newtheorem{thm}[equation]{Theorem}

\newtheorem{lem}[equation]{Lemma}
\newtheorem{cor}[equation]{Corollary}

\newtheorem{prop}[equation]{Proposition}

\newtheorem{hyp}[equation]{Hypotheses}

\newtheorem{claim}[equation]{Claim}

\newtheorem{rmk}[equation]{Remark}

\theoremstyle{definition}
\newtheorem{defn}[equation]{Definition}

\newcommand{\cX}{{\check{X}}}

\newcommand{\frr}{{\mathfrak r}}

\newcommand{\bbS}{{\mathbb S}}

\newcommand{\CO}{{\mathcal{O}}}
\newcommand{\CV}{{\mathcal{V}}}

\newcommand{\CJ}{{\mathcal{J}}}
\newcommand{\CF}{{\mathcal{F}}}
\newcommand{\CE}{{\mathcal{E}}}
\newcommand{\cL}{{\mathcal{L}}}
\newcommand{\CW}{{\mathcal{W}}}

\newcommand{\SK}{{\bbS_{K}}}
\newcommand{\SKt}{{\bbS_{K}^{\rm tor}}}

\newcommand{\ZZ}{{\mathbb Z}}
\newcommand{\TT}{{\mathbb T}}
\newcommand{\Gm}{{\mathbb G}_m}

\newcommand{\CB}{{\mathcal{B}}}
\newcommand{\Qpb}{{\bar{\QQ}_p}}

\newcommand{\CS}{{\mathcal{S}}}

\newcommand{\CT}{{\mathcal{T}}}
\newcommand{\CG}{{\mathcal{G}}}

\newcommand{\G}{{\Gamma}}

\newcommand{\fm}{{\mathfrak{m}}}

\newcommand{\ra}{{~\rightarrow~}}

\newcommand{\CH}{{\mathcal{H}}}

\newcommand{\QQ}{{\mathbb Q}}

\newcommand{\Qbar}{{\overline{\mathbb Q}}}

\newcommand{\FF}{{\mathbb F}}

\newcommand{\RR}{{\mathbb R}}
\newcommand{\ad}{{\mathbf A}}
\newcommand{\af}{{{\mathbf A}_f}}
\newcommand{\afp}{{{\mathbf A}_f^p}}
\newcommand{\afS}{{{\mathbf A}_f^S}}
\newcommand{\CC}{{\mathbb C}}
\newcommand{\MM}{{\mathbb M}}

\newcommand{\Qp}{{{\mathbb Q}_p}}
\newcommand{\Zp}{{{\mathbb Z}_p}}
\newcommand{\Qq}{{{\mathbb Q}_q}}

\newcommand{\Hom}{{\mathrm{Hom}}}
\newcommand{\adj}{{\mathrm{ad}}}
\newcommand{\Aut}{{\mathrm{Aut}}}
\newcommand{\gr}{{\mathrm{gr}}}
\newcommand{\Endo}{{\mathrm{End}}}
\newcommand{\isoarrow}{{~\overset\sim\longrightarrow}}
\newcommand{\br}{\bar{r}}

\newcommand{\Frob}{\mathrm{Frob}}
\newcommand{\Lie}{\mathop{\rm Lie}}
\newcommand{\Spec}{{\mathrm{Spec}}}
\newcommand{\Gal}{{\mathrm{Gal}}}

\begin{document}

\title[The Taylor-Wiles method for coherent cohomology, II] % (optional, use only with long paper titles)
{The Taylor-Wiles method for coherent cohomology, II}

\author[Stanislav Atanasov]{Stanislav Atanasov}

\author[Michael Harris]{Michael Harris}
%\address{Institut de Math\'ematiques de
%Jussieu, U.M.R. 7586 du CNRS; UFR de Math\'ematiques \\ Universit\'e Paris-Diderot Paris 7}

\smallskip

%\end{classification}

%\keywords{Galois representation, Shimura variety, special values of $L$-functions}
%\end{keywords}
\subjclass[2010]{}
\thanks{The research leading to these results has received funding from the European Research Council under the European Community's Seventh Framework Programme (FP7/2007-2013) / ERC Grant agreement no. 290766 (AAMOT).  The second author was partially supported by NSF grants DMS-1404769 and DMS-2001369}

\maketitle

\setcounter{tocdepth}{1}
\tableofcontents

\begin{abstract}
We show that, under certain specific hypotheses, the Taylor-Wiles method can be applied to the  cohomology of a Shimura variety $S$ of PEL type attached to a unitary similitude group $G$, with coefficients in the coherent sheaf attached to an automorphic vector bundle $\CF$ , when $S$ has a smooth model over a $p$-adic integer ring.  As an application, we show that, when the hypotheses are satisfied, the congruence ideal attached to a coherent cohomological realization of an automorphic Galois representation is independent of the signatures of the hermitian form to which $G$ is attached.  We also show that the Gorenstein hypothesis used to construct $p$-adic $L$-functions  in \cite{EHLS}, as elements of Hida's ordinary Hecke algebra, is valid rather generally.

The present paper generalizes the main results of the article \cite{H13}, which treated the case when $S$ is compact.    As in the previous article, the starting point is a theorem of Lan and Suh that proves the vanishing of torsion in the cohomology under certain conditions on the parameters of the bundle $\CF$  and the prime $p$.   Most of the additional difficulty in the non-compact case
is related to showing that the contributions of boundary cohomology are all of Eisenstein type.  We also need to show that the coverings giving rise to
the diamond operators can be extended to  \'etale coverings of appropriate toroidal compactifications.  

\end{abstract}

\section*{Introduction}

The Taylor-Wiles method was introduced as a way of proving automorphy lifting theorems, and it has been applied in a variety of situations to show, under appropriate technical hypotheses, that all liftings to characteristic zero of an $n$-dimensional Galois representation $\bar{\rho}$ with coefficients in a finite field are attached to automorphic forms, provided one knows that one of the liftings is automorphic.  
As a by-product of the method, one finds -- in the {\it minimal} case -- that the localization of the relevant Hecke algebra at the maximal ideal attached to $\bar{\rho}$ is a local complete intersection, and that the module of automorphic forms is free over the localized Hecke algebra $\TT$.   These properties of the Hecke algebra and its module of automorphic forms are of independent interest.  They are used, for example, in the construction of $p$-adic $L$-functions in \cite{EHLS}.  

The primary aim of the present paper, like its predecessor \cite{H13}, is to prove that these properties are satisfied, in some generality, by the {\it coherent} cohomology of the Shimura varieties attached to unitary similitude groups.  To this end we develop the Taylor-Wiles method in this situation
and indicate some of the applications of the local complete intersection property to automorphic forms.  

The reader should understand that, in contrast to most papers on the Taylor-Wiles method, our aim is {\it not} to show that certain Galois representations are automorphic (or potentially automorphic).   This is already known for the kinds of Galois representations that arise in the cohomology of the PEL Shimura varieties that are the subject of these papers.  In contrast, the primary purpose of the present paper, as well as its predecessor \cite{H13}, is to obtain results about the structure of the localized Hecke algebra $\TT$, with a view to applications to arithmetic problems.  Thus the main results of the paper are those contained in \S 4 and \S 5.  Thus, let $F$ be a CM field in which every prime dividing $p$ splits over the maximal totally real subfield $F^+$ of $F$, and let $\pi$ be a cuspidal cohomological automorphic representation of $GL(n)_F$.  We assume the contragredient of $\pi$ is isomorphic to  $\pi^c$, where $c$ is the non-trivial element of $Gal(F/F^+)$.  When $V$ is an $n$-dimensional hermitian space over $F$, and $G$ its unitary similitude group, we consider the $L$-packet $\Pi_V(\pi)$ of automorphic representations of discrete series type that base change to $\pi$, and we assume $\Pi_V(\pi)$ is non-empty; the conjugate duality hypothesis is a necessary condition for this, but there may also be local conditions at places where $G$ is not quasi-split (see \cite{L}).  We temporarily denote by $S_V$  the corresponding Shimura variety.  We specifically prove -- under the hypotheses that allow us to apply the Taylor-Wiles method -- that
\begin{itemize}
\item The exterior power relation between the $\Qbar_p$-valued cohomology spaces of $S_V$, for $V$ of different signatures, attached to a given cuspidal automorphic  $\pi$, is in fact integral over the Hecke algebra localized at the maximal ideal corresponding to $\pi$ (Proposition \ref{ext});
\item The congruence ideals for the localized Hecke algebras attached to $\Pi_V(\pi)$ coincide for different $V$ (Proposition \ref{congruenceV});
\item If $\pi$ is ordinary at $p$ and $\pi_V^0 \in \Pi_V(\pi)$ is the holomorphic member, then the Gorenstein hypothesis of \cite{EHLS} is valid for the Hida family through $\pi_V^0$; in particular, the $p$-adic $L$-function attached to this Hida family in \cite{EHLS} takes values in the localized (big) Hecke algebra.
\end{itemize}
These results are all obtained as consequences of the Taylor-Wiles method for coherent cohomology.

In the article \cite{H13}, the vanishing theorem \cite{LS1} of Lan and Suh was used to show that the Taylor-Wiles method can be applied to the (coherent) cohomology of a compact Shimura variety $S$, attached to a reductive algebraic group $G/\QQ$, with coefficients in an automorphic vector bundle $\CF$.  Assuming 
\begin{itemize}
\item[(i)] the parameter of
$\CF$ is in the range to which the vanishing theorem applies; 
\item[(ii)]  $S$ has a smooth model $\SK$ over a $p$-adic integer ring $\CO$; here the subscript $K$ denotes a level subgroup, and 
\item[(iii)] there is a theory of Galois representations attached to cohomological automorphic forms on $S$,  
\end{itemize}
the cohomology $H^{\bullet}(\SK,\CF)$, with coefficients in $\CO$  is concentrated in a single degree  $q$ and $H^{q}(\SK,\CF)$, localized at a maximal ideal as above, is free over the localized Hecke algebra.    
The crucial input in the Taylor-Wiles method, applied to a reductive group $G$, is the action of a product of groups of the form 
$$\Delta_{Q_N} = (\prod_{q \in Q_N} K_0(q)/K_1(q))_p$$
where $Q_N$ is a finite set of {\it Taylor-Wiles primes} of the base field $F$ (usually totally real or a CM field), with $Nq \equiv 1 \pmod{p}^N$, $K_0(q)$ and $K_1(q)$ are open compact subgroups of $G(F_q)$\footnote{This doesn't literally make sense unless $F = \QQ$, but we ignore this in the present exposition}
that generalize the classical congruence subgroups $\Gamma_0(q)$ and $\Gamma_1(q)$ respectively, so that $K_0(q)/K_1(q)$ is isomorphic to the multiplicative group of the residue field $k(q)$, and the subscript $p$ denotes the maximal quotient of exponent $p^N$.    Let $C^{Q_N}$ denote the kernel 
of the map $\prod_{q \in Q_N} K_0(q)/K_1(q) \rightarrow \Delta_{Q_N}$, and let $\SK(Q_N)$ denote the quotient by $C^{Q_N}$ of the Shimura variety at level 
$K(Q)$, which is the subgroup obtained by replacing the local component of $K$ at a prime  $q \in Q$ by $K_1(q)$.  Provided the original level $K$ is sufficiently deep, the action of $\Delta_{Q_N}$ on $\SK(Q_N)$ is {\it free}.  A theorem of  S. Nakajima \cite{N}, combined with the Lan-Suh vanishing theorem, then implies that $H^{q}(\SK(Q_N),\CF)$ is a free $\CO[\Delta_{Q_N}]$-module.  This, together with item (iii) above, allows us to prove (under favorable hypotheses on $\bar{\rho}$) that the original $H^{q}(\SK,\CF)$ is free over the localized Hecke algebra $\TT$.

The present paper is a sequel to \cite{H13}; its purpose is to extend these results to non-compact Shimura varieties, specifically those attached to unitary similitude groups of hermitian vector spaces over CM fields.   To this end we apply the vanishing theorem of a different Lan-Suh paper \cite{LS2}, which extends their previous results to non-compact Shimura varieties.  

The project faces two immediate difficulties.  In the first place, coherent cohomology in the non-compact case is computed on toroidal compactifications, not on the open Shimura variety.  In order to apply Nakajima's result, we need to show that the action of the group denoted $\Delta_{Q_N}$ above on some toroidal compactification 
$\SK(Q_N)^{tor}$ of $\SK(Q_N)$ is still free; in other words, the $\SK(Q_N)^{tor}$ is an \'etale cover of its quotient by $\Delta_{Q_N}$.  This is certainly not true for an arbitrary toroidal compactification, and one of the main objectives of the present paper is to construct toroidal compactifications of the integral models constructed in \cite{Lan} on which $\Delta_{Q_N}$ acts freely.   In the classical case of elliptic modular curves, the fact that the complete modular curve $X_1(q)$ is (usually) \'etale over $X_0(q)$ is a well known special feature of this particular pair of congruence subgroups.   The constructions in the present paper make extensive use of the results of \cite{Lan} and \cite{Lan2}.

The second difficulty is that the Lan-Suh vanishing theorem in the non-compact case applies only to {\it interior coherent cohomology}, in other words the image $H^{\bullet}_!(\SK,\CF)$ of the cohomology of $\SKt$ with coefficients in the subcanonical extension $\CF^{sub}$ in the cohomology of the canonical extension 
$\CF^{can}$.  Thus, although Nakajima's result applies to the canonical extension, one cannot apply it directly to the problem
at hand.   The freeness of the action allows us to represent $H^{\bullet}(\SKt,\CF^{can})$ as the cohomology of a perfect complex 
$R\Gamma^{\bullet}$ of $\CO[\Delta_{Q_N}]$-modules with an action of the (non-localized) Hecke algebra $T$, which we define, following Khare and Thorne \cite{KT}, as an algebra of endomorphisms of the object $R\Gamma^{\bullet}$ in an appropriate derived category.  Thus, even though our final results concern a Hecke algebra module that occurs in a single degree of cohomology, we find ourselves obliged to work, at least in the first stages, with the extension of the Taylor-Wiles method to cohomological complexes introduced by Calegari and Geraghty \cite{CG}, and pursued by Hansen \cite{Han} as well as Khare-Thorne.

We then localize $T$ at a non-Eisenstein maximal ideal $\fm$, and write $\TT = T_{\fm}$.  Using the results of \cite{HLTT} and the fact, due to Scholze, Boxer, Pilloni-Stroh, and Goldring-Koskivirta \cite{Sch,Box,PS,GK}, that there is a theory of Galois representations attached to {\it torsion} cohomology classes as well, we conclude that the localization at $\fm$ of $H^{\bullet}(\SKt,\CF^{can})$ coincides with the localization of $H^{\bullet}_!(\SK,\CF)$.  Thus it is concentrated in a single degree $q$, and is $p$-torsion free, when the Lan-Suh results apply.  It then follows easily that the localized cohomology is free over $\CO[\Delta_{Q_N}]$, and we conclude the argument as in \cite{H13}.  

This material occupies the first three sections of the paper.  Sections \ref{modules} and \ref{ordinarymf} develop applications of these results. 
In Section \ref{modules}, we apply the results for coherent cohomology in order to prove that the interior de Rham and $p$-adic \'etale cohomology are also free over the Hecke algebra, after localization at a non-Eisenstein prime, and again in the range to which the results of \cite{LS2} apply.  As $V$ varies among hermitian spaces of a fixed dimension $n$, the Galois representations on the corresponding $p$-adic \'etale cohomology groups are related by a formula conjectured by Langlands that has essentially been proved \cite{ScholzeShin}  using the methods of Kottwitz and the stable trace formula.  Roughly speaking, up to twist by an explicit character, all of the cohomology groups can be obtained from a few of them by tensor operations.  We show that the analogous relations hold over the localized Hecke algebra.  We also observe that the congruence ideal for the automorphic representations realized
on the unitary Shimura varieties attached to varying $V$ depends only on the associated Galois representation and not on $V$.  These results are only proved, of course, when the Taylor-Wiles hypotheses are valid and assuming the parameters are in the range of the Lan-Suh vanishing theorem.

In section \ref{ordinarymf} we  specialize to degree $0$ but apply the results to holomorphic modular forms varying in ordinary Hida families.  We prove that, 
when the Taylor-Wiles hypotheses apply to residual representation of a Hida family, then the family satisfies the Gorenstein hypotheses 
used in \cite{EHLS} to construct $p$-adic $L$-functions as elements of Hecke algebras.   

The last two sections are devoted to the proof that the contribution of the boundary cohomology is purely Eisenstein.   This roughly follows the pattern of
the analysis of boundary coherent cohomology in \cite{HZ94,HZ94b,HZ01}.   However, some of the arguments in those papers are analytic in nature, or are only valid for cohomology with characteristic zero coefficients, and there is no reason to suppose that the results are as clean as in the analytic setting.  Fortunately, for our purposes we can be satisfied with qualitative results.  To prove these we  make extensive use of Lan's study of the geometry of integral models of Kuga varieties in \cite{Lan,Lan2,Lan3}.  We also borrow some ideas of Newton and Thorne from \cite{NT}.

An appendix explains the translation between the boundary cohomology of connected Shimura varieties, as developed in \cite{HZ94,HZ94b}, which is more convenient for the local study of the boundary, and the adelic theory of \cite{HZ01}, which simplifies level-raising arguments.

\subsection*{Conventions}
Let $F$ be a CM field, quadratic over a totally real field $F^+$ and contained in a fixed algebraic closure $\Qbar$ of $\QQ$, $d = [F^+:\QQ]$  We let
$\Gamma_F = \Gal(\Qbar/F)$, $\Gamma_{F^+} = \Gal(\Qbar/F^+)$.  
Let $S_\infty$, resp. $S^+_{\infty}$ denote the set of complex embeddings of $F$, respectively real places
of $F^+$, and let $\tilde{S}_{\infty}$ denote a CM type for $F$, i.e. a choice for each $v \in S^+_{\infty}$ of
a complex place $\tilde{v} \in S_\infty$.    Let $S(F^+)$ denote the set of primes of $F^+$ that ramify in $F$.

Our modular forms and deformation problems will be defined over a $p$-adic integer ring $\CO$, with maximal ideal
$\fm_\CO$.  It will always be taken
big enough to contain the rings of definition of all local deformation problems.

\section{Shimura varieties attached to unitary groups}

\subsection{Notation for automorphic vector bundles}\label{avbnotation}

Let $V$ be an $n$-dimensional
space over $F$ with nondegenerate hermitian form $(\bullet,\bullet)$,  let $U = U(V)$ be its unitary group, and define the reductive group
$G$ over $\QQ$ by its values on $\QQ$-algebras $R$:
\begin{equation}\label{similitude} G(R) = \{g \in GL(V\otimes_{\QQ}R) ~|~ (g(u_1),g(u_2)) = \nu(g)(u_1,u_2) \}\end{equation}
for all $u_1, u_2 \in V\otimes_{\QQ}R$ for some $\nu(g) \in R^{\times}$.
As in \cite{H19} we extend $G$ to a Shimura datum $(G,X)$ so that 
$$S(G,X)(\CC) := \varprojlim_{K \subset G(\af)}  G(\QQ)\backslash X \times G(\af)/K = \varprojlim_{K \subset G(\af)} S_K(G,X)(\CC),$$
where $K$ varies over open compact subgroups of $G(\af)$, is the set of complex points of a Shimura
variety $S(G,X) = \varprojlim S_K(G,X)$ of PEL type, with reflex field $E = E(G,X)$.  We let $d_V = \dim S_K(G,X)$.  
The CM type $\tilde{S}_{\infty}$ is part of the data defining $X$.  For each $v \in S^+_{\infty}$,
there is a partition $r_v + s_v = n$ such that the hermitian space $V\otimes_{F,\tilde{v}}\CC$ has signature $(r_v,s_v)$; then
$U(F^+_v) \isoarrow U(r_v,s_v)$, with the ambiguity between $U(r_v,s_v)$ and $U(s_v,r_v)$ resolved by
the choice of $\tilde{v}$.

In what follows we always assume $K = \prod_q K_q$, where $q$ runs over rational prime numbers and $K_q \subset G(\Qq)$ is open compact.  We fix an odd prime $p$ for the remainder of this paper and assume $p$ is not in the set $S = S(K)$ of bad primes for $K$, defined in \S \ref{typedata} below; this can be weakened slightly, but we do not bother to treat the more general situation.

For any $q$, we break up the set of primes $D(q)$ of $F^+$ dividing $q$ into the subsets $D_+(q) \coprod D_-(q)$, where $v \in D_+(q)$ if and only if $v$ is split in $F/F^+$.  For each $v \in D(q)$ choose a prime $w(v)$ of $F$ dividing $v$.  Then there is a natural isomorphism
\begin{equation}\label{divisors} G_q := G(\Qq) \isoarrow \prod_{v \in D_+(q)} GL(n,F_{w(v)}) \times G_{-,q} \end{equation}
where we set $V_- = \prod_{v \in D_-(q)} V\otimes_{\QQ}F^+_v$ and
$$G_{-,q} = \{g \in GL(V_-) ~|~ (g(u_1),g(u_2)) = \nu(g)(u_1,u_2) \text{ for some } \nu(g) \in \Qq^{\times}\},$$
where $u_1$, $u_2$ vary over vectors in $V_-$.  For all $q$ we assume $K_q$ admits a factorization
\begin{equation}\label{factorkq}  K_q = \prod_{v \in D_+(q)} K_v \times K_{-,q} \end{equation}
where $K_v \subset GL(n,F_{w(v)})$ and $K_{-,q} \subset G_{-,q}$ are open compact subgroups.

\subsubsection{Type data}\label{typedata} It was pointed out in \cite[Remark 6.9]{H13} that the hypotheses of that paper restricted ramification of
the residual representation $\bar{\rho}$ at places of $F^+$ not dividing $p$.  In this paper those restrictions are somewhat relaxed, although we only consider
deformation problems that involve minimal ramification.
Let $S = S(K)$ be the set of primes $q$ where either $K_v$, for some $v \in D_+(q)$, or $K_{-,q}$, is not hyperspecial maximal; we call $S$ the set of {\it ramified primes} for $K$.  The set $S$ contains the set $S(G)$ of primes $q$ at which the group $G$ is ramified, and in particular contains the rational primes divisible by (the set) $S(F^+)$ of primes of $F^+$ that ramify in $F/F^+$.  
At primes $q \in S(G)$ we take $K_{-,q}$ to be special maximal. 
If $q \in S \setminus S(G)$ we assume $K_{-,q}$ is a hyperspecial maximal subgroup of $G_{-,q}$.  For $v \in D_+(q)$, whether or not
$q \in S(G)$, we will choose a set of {\it type data}, namely quadruples $(K_v^+,K_v,\Yu_v,\Lambda_{\Yu_v})$, with $K_v \subset K^+_v$ a pair of open compact subgroups, $\Lambda_{\Yu_v}$
a finite free $\CO$-module and
\begin{equation}\label{BKtype}
\Yu_v:  K_v^+/K_v \ra \Aut(\Lambda_{\Yu_v})
\end{equation}
a representation that determines the inertial type of the corresponding representation of $GL(n,F^+_v)$, as in \cite{BK99}.

We will also use $S$ to denote the set of $v$ dividing primes $q \in S(K)$ where the factor $K_v$ of $K_q$ is not hyperspecial maximal; this in particular allows us to use the same notation even when $q$ is ramified in $F^+$.

Our level subgroups $K$ will always be assumed to be {\it neat}, as in \cite{H89,P}, so that  the Shimura variety $S_K(G,X)$ is smooth and we can choose toroidal compactifications that are smooth and projective.  The set $S$ will be assumed to contain a prime $\frr$ all of whose divisors in $F^+$ split in $F/F^+$.   We will choose one such divisor $\frr_0$ so that $K_{\frr_0}$ is the subgroup of a standard maximal compact subgroup  of $GL(n,F^+_{\frr_0})$ 
defined by
\begin{equation}\label{neat0} K_{\frr_0} = \{k \in GL(n,\CO_{\frr_0})~|~ k \equiv u \pmod \fm_{\frr_0}\},
\end{equation}
where $\fm_{\frr_0}$ is the maximal ideal in the integer ring and $u$ is upper triangular unipotent.  This suffices to guarantee that $K$ is neat, and $\frr_0$ will be chosen
so that the deformation problems to be studied below are minimal and unrestricted at $\frr$.  In the statement of \cite[Theorem 6.8]{H13} this set is called $S_a$.\footnote{Erratum for \cite{H13}:  In the published version of \cite{H13}  the set $S_a$ was mentioned appropriately in the statement of its main theorem but its
definition was inadvertently omitted. As already mentioned, this condition is required in order to guarantee smoothness.  Theorem \ref{TWtheorem} requires as well that the
deformation condition at such an $\frr$ be unrestricted and minimal.  We choose $\frr$  as on \cite[p. 111]{CHT} (where it is called $S_a$) or as on
\cite[p. 893]{Th} (where it is called $v_1$).  The hypothesis on $\frr$ -- that $[F(\zeta_\frr):F] > n$ -- guarantees that $\frr$ satisfies condition (1)
of Theorem \ref{TWtheorem}.}

We recall the standard constructions of automorphic vector bundles.  The points $x \in X$ are identified with homomorphisms 
$$h_x: \mathbb{S} = R_{\CC/\RR} \mathbb{G}_{m,\CC} \ra G_{\RR}$$
satisfying Deligne's list of axioms [De, 2.1.1].   The centralizer in $G(\RR)$ of $h_x$ is a reductive
group $K_x$ whose intersection with the derived subgroup $G^{der}(\RR)$ of $G(\RR)$
is a maximal compact connected subgroup.  Deligne's axiom (2.1.1.1) concerns the adjoint action
$$Ad\circ h_x:  \mathbb{S} \ra GL(\Lie(G)_{\CC})$$
that yields an eigenspace decomposition (the {\it Harish-Chandra decomposition})
\begin{equation} \Lie(G)_{\CC} = \mathfrak{p}_x^- \oplus \Lie(K_x)_{\CC} \oplus \mathfrak{p}_x^+. \end{equation}
Here $z \in \mathbb{S}(\RR) \isoarrow \CC^{\times}$ acts trivially on $\Lie(K_x)$ and as
$(z/\bar{z})$ (resp. $\bar{z}/z$) on $\mathfrak{p}_x^+$ (resp. $\mathfrak{p}_x^-$).
%%$\mathfrak{p}_x^- = \Lie(G)^{1,-1}$
The Lie subalgebras $\mathfrak{p}_x^-$ and $\mathfrak{p}_x^+$ are naturally identified, respectively, with the
anti-holomorphic and holomorphic tangent spaces of $X$ at $x$. 

Let $\check{X}$ denote the compact dual of $X$, and $X \hookrightarrow \check{X}$ the Borel
embedding.  Concretely,  $\check{X}$ is a flag variety of maximal parabolic subgroups
of $G_{\CC}$ and the image in $\check{X}$ of $x \in X$ is a maximal parabolic $P_x$ with Levi
subgroup $K_x$.  
In [H85] it is explained how to define a canonical $E = E(G,X)$-rational structure on the flag variety $\cX$, following Deligne,
and how to define a functor $\CV \mapsto [\CV]$ from $G$-equivariant vector bundles on $X$ to $G(\afp)$-equivariant
vector bundles on $S_{K}(G,X)$.   The latter are 
called {\it automorphic vector bundles}.   The functor is compatible with the $E$-structure in the sense that, if 
$\CV$ is defined as $G$-equivariant vector bundle over a field $E(\CV)$ (which can always be taken to be
a number field), then for any $\sigma \in \Gal(\bar{E}/E)$, we have
\begin{equation} \sigma[\CV] = [\sigma(\CV)]. \end{equation}

The bracket notation in the previous paragraph was introduced in order to make reference to the functor.
Automorphic vector bundles will in general be denoted $\CF$, and we let $E(\CF)$ denote a field
of definition for the corresponding equivariant vector bundle over $\check{X}$.

Fix a point $x \in X$ with stabilizer $P_x \subset G_{\CC}$; thus
 \begin{equation} \Lie(P_x) = \Lie(K_x)_{\CC} \oplus \mathfrak{p}_x^+ \end{equation}
in the Harish-Chandra decomposition.
There is a natural equivalence
of categories between $G$-equivariant
vector bundles $\CV$ on $X$ and finite-dimensional representations $(\tau,W_{\tau})$ of $P_x$;  $W_{\tau}$
is the fiber of $\CV$ at $x$, and $\tau$ is the isotropy representation.  
 In particular, the natural
representation of $P_x$ 
$$ad^+:  P_x \ra K_x \ra \Aut(\mathfrak{p}_x^+),$$
where the second arrow is the adjoint representation, defines an automorphic vector bundle canonically
isomorphic to the tangent bundle $\CT_{S(G,X)}$.  Likewise $\wedge^{top}(ad^+)^{\vee}$, the dual of
the top exterior power of the adjoint action on $\mathfrak{p}_x^+$, defines the canonical bundle $\Omega^{top}_{S(G,X)}$
as an automorphic vector bundle.   

A representation $(\tau,W_{\tau})$ of
$K_x$ extends trivially to a representation of $P_x$ and thus defines an automorphic vector bundle
$\CF = \CF_{\tau}$  on $S_{K_p}(G,X)$ whose fiber at a point beneath $x \times g$ for any $g \in G(\af)$ 
can be identified with $W_{\tau}$.  The automorphic vector bundle $\CF$ can also be identified
with the family of bundles $\CF_K$ on $S_K(G,X)$.   
If the family $\CF_K$ extends to a $G(\afp)$-equivariant family of vector bundles on a family of $\CO$-integral models
$\bbS_K$ of $S_K(G,X)$, where $\CO$ is some $p$-adic integer ring, we denote the extension $\CF_{K,\bbS}$.  In the applications $\bbS_K$ is a moduli space for abelian varieties with PEL structure and natural extensions will be specified in terms of this identification.

\subsubsection{Notation for highest weights}  Notation is as in \cite[\S 2]{H13}.  We choose a maximal torus $T \subset K_x$; 
then $T$ is also a maximal torus in $G$.   
Let $\Phi = \Phi(G,T)$ denote the set of roots of $G$ relative to $T$, and let $\Phi^+$ be a system of positive
roots compatible with $P_x$, i.e. containing the roots of $\mathfrak{p}_x^-$; let $X^+(T)$ (resp. $X_x^+(T)$) denote the set of dominant weights for $G$ (resp. for $K_x$) relative to this choice.   Let $W = W(G,T)$ be the
absolute Weyl group and let $W^x \subset W$ be the set of Kostant representatives.    

\begin{defn}\label{psmall}  (cf. \cite[2.32]{LS1}) Let $\mu \in X^+(T)$.     Say $\mu$ is {\it $p$-small}, resp.  {\it $p$-small relative to $x$},
if
$$\langle \mu + \rho,\alpha^{\vee}\rangle ~~\leq ~~p, \forall ~\alpha \in X^+(T) \text{ (resp.  $\forall ~\alpha \in X^+_x(T)$)} $$
and if $p > |\mu|_L$, where if we write $\mu = (\mu_\sigma, \sigma \in \tilde{S}_\infty)$,
$$|\mu|_L = \sum_{\sigma \in \tilde{S}_\infty} |\mu_\sigma|$$
 in the notation of \cite[Definition 3.2]{LS1}. (In \cite{LS1} this is called ``$p$-small for the geometric realization of Weyl's construction.")
\end{defn}

%%%%%%%%%%%%%%%%%%%%%%%%%%%%%%%%%%%%%%%%%%%%%%%%%%%%%%%%%%%%%%%%%%%%%%%%%%%%%%%%%%%%%%%%
%%%%%%%%%%%%%%%%%%%%%%%%%%%%%%%%%%%%%%%%%%%%%%%%%%%%%%%%%%%%%%%%%%%%%%%%%%%%%%%%%%%%%%%%

\subsection{Integral models}  
\label{sec:integral-models}
The Shimura datum $(G,X)$ is of PEL type and $G$ is of type $A$, so by \cite{Ko} we may consider integral models of the Shimura variety $S_K(G,X)$ over $p$-adic integer rings, at least when $p$ is unramified for $G$ and $K$. Let $g=n\cdot [F^+:\QQ]$.   

\begin{prop}
\label{prop:integral-model-open}
Let $p$ be a prime at which $G$ is unramified, and fix a hyperspecial maximal compact $K_p\subseteq G(\QQ_p).$ Then the Shimura variety $S_K(G,X)$ admits a smooth model $\bbS_K$ over a $p$-adic integer ring $\CO$ for all neat $K\supset K_p.$ More precisely, if $K$ is neat and contains $K_p,$ then up to replacing $K$ by a normal subgroup $K'$ of finite index, for any prime $v$ in $E(G,X)$ above $p$ there exists a smooth moduli scheme $\bbS_{K'}$ over $\Spec(\CO_v)$ of abelian varieties of dimension $g$, with PEL structure defined in terms of the hermitian space $V$, whose generic fiber is isomorphic to $S_{K'}(G,X)$. The quotient $\bbS_{K'}$ by $K/K'$ supplies the integral model $\bbS_K.$
\end{prop}

\begin{rmk}\label{Hasse}  It is well known that $S_K(G,X)$ is a subvariety of the moduli space parametrizing all quadruples $(A,i,\lambda,\kappa)$ where
$A$ is an abelian variety of dimension $g$ with endomorphisms by (a subring of) $F$ determined by $i:  F \hookrightarrow \Endo(A)$, polarization $\lambda$,
and level structure $\kappa$, all satisfying the usual hypotheses adapted to the signatures of $V$ at archimedean places of $F$.  We denote this moduli space
$M_K(G,X)$; it is a union of a finite number of copies of $S_K(G,X)$, corresponding to the number of global forms of $G$ that are locally isomorphic everywhere.
The moduli space $M_K(G,X)$ plays no separate role in the theory developed in this paper.
\end{rmk}

%%%%%%%%%%%%%%%%%%%%%%%%%%%%%%%%%%%%%%%%%%%%%%%%%%%%%%%%%%%%%%%%%%%%%%%%%%%%%%%%%%%%%%%%
%%%%%%%%%%%%%%%%%%%%%%%%%%%%%%%%%%%%%%%%%%%%%%%%%%%%%%%%%%%%%%%%%%%%%%%%%%%%%%%%%%%%%%%%
\subsection{Terminology for toroidal compactifications}\label{termin}

For any level subgroup $K$, the Shimura variety $S_K(G,X)$ has a family of toroidal compactifications, each one attached to a collection $\Sigma$ of combinatorial data adapted to $K$.  The precise definition of $\Sigma$ is recalled below. The corresponding toroidal compactification is denoted $S_K(G,X)_{\Sigma}$, or $S_K(G,X)^{tor}$ when we don't need to specify $\Sigma$.  When $K$ is allowed to vary through a collection $\CB$ of open compact subgroups of $G(\af)$ , we choose $\Sigma(K)$ adapted to $K \in \CB$ in such a way that, if $K' \subset K$ with both $K'$ and $K$ in $\CB$, the natural covering map $S_{K'}(G,X) \ra S_{K}(G,X)$ extends to a map $$S_{K'}(G,X)_{\Sigma(K')} \ra S_{K}(G,X)_{\Sigma(K)}$$ of toroidal compactifications. 

The combinatorial data $\Sigma=\cup \Sigma_F$ adapted to the neat level $K$ is a collection of fans, $\Sigma_F$, one for each rational boundary component $F$.\footnote{We will make an effort to use the letter $F$ for boundary components and for number fields in different paragraphs, so there should be no confusion.}
 Each $F$ corresponds to its stabilizer $P_F$, which is a maximal rational parabolic inside $G.$ The fan $\Sigma_F$ gives a polyhedral cone decomposition of a
 partial compactification $\bar{C}_F$ of a certain cone $C_F$ inside $U_F(\RR)$ that is open, convex, and self-adjoint with respect to a $\QQ$-rational positive definite quadratic form, where $U_F$ is the center of the unipotent radical of $P_F.$   We examine the $\Sigma_F$ and the inclusion
 $C_F \subset \bar{C}_F$ more closely below.

If $\Sigma$ is $K\cap G(\QQ)$-admissible in the sense of \cite[Definition 5.1]{AMRT}, the compactification $S_K(G,X)_\Sigma$ is smooth, and if $\Sigma$ is furthermore defined by cocores then $S_K(G,X)_\Sigma$ is projective by Tai's theorem \cite[IV, \S2]{AMRT}. Lastly, we assume that $\Sigma$ is such that 
$$\partial S_K(G,X)_\Sigma := S_K(G,X)_\Sigma - S_K(G,X)$$ 
is a divisor with normal crossings. This is equivalent to the hypothesis that for all $\sigma\in\Sigma$, the semigroup $\sigma\cap (U_F(\QQ)\cap K)$ is generated by a subset of a basis for the free abelian group $U_F(\QQ)\cap K$.  We will always choose  $\Sigma$ that satisfy all of the above conditions; such $\Sigma$ exist and are constructed for instance in \cite{H89}.

Unlike the anisotropic case treated in \cite{H13}, to apply the vanishing theorem of Lan-Suh for non-compact Shimura varieties, we need work with a smooth proper integral model $\bbS_K^{tor}$, over a $p$-adic integer ring $\CO$, of the toroidal compactification, for which $\bbS_K(G,X)$ of \S\ref{sec:integral-models} embeds as an open dense subscheme. In the case when $G$ is a symplectic similitude group this has been constructed in the book \cite{FC} by Faltings-Chai and subsequently extended to all PEL type Shimura varieties by K.-W. Lan in his thesis \cite{Lan}. The comparison between the algebraic (moduli) construction of these compactification and the analytical one presented above is established in \cite{Lan12}. Below is an abridged summary of the relevant results.

\begin{thm}(Lan)\label{toroidalintegral}
Let $p$ be a prime at which $G$ is unramified, and fix a hyperspecial maximal compact $K_p\subseteq G(\QQ_p).$ Then for all neat $K\supset K_p$, there is a compatible choice of admissible smooth rational polyhedral data $\Sigma$ (see \cite[Def. 6.3.3.2]{Lan}) such that the $\bbS_K(G,X)_\Sigma$ is a smooth proper scheme over a $p$-adic integer ring $\CO$ which contains $\bbS_K(G,X)$ of Proposition~\ref{prop:integral-model-open} as an open dense subscheme. Furthermore, $\Sigma$ may be chosen so that $\bbS_K(G,X)_\Sigma - \bbS_K(G,X)$, viewed as a closed reduced subscheme, is a divisor with normal crossings. There is a canonical strata-preserving isomorphism between the basechange to $\Spec(\CC)$ of this integral model and the classical analytical construction as in \cite{AMRT}.
\end{thm}

\subsubsection{Parabolic strata}

The boundary divisor $\partial S_K(G,X)_\Sigma$ has a closed covering indexed by maximal standard rational parabolic subgroups of $G$, defined
as follows.  Let $\bbS_K^{min}$ denote the minimal compactification of $\bbS_K$ over $Spec(\CO)$, as in \cite[Theorem 7.4.2.1]{Lan}.  
The boundary $\bbS_K^{min} \setminus \bbS_K$ decomposes as a disjoint union of locally closed strata $\partial^P\bbS_K^{min}$ indexed by standard maximal rational parabolic subgroups $P \subset G$.  Each $\partial^P\bbS_K^{min}$ is a union of certain strata
denoted $\mathbf{Z}_{[(\Phi_{\mathcal{H}},\delta_{\mathcal{H}})]}$ in Lan, where $P$ is the stabilizer of the filtration on the hermitian vector space $V$ concealed in the notation $(\Phi_{\mathcal{H}},\delta_{\mathcal{H}})$ (see \cite[\S 5.4]{Lan} for an explanation).

There is a canonical morphism  \cite[Theorem 7.4.2.1, 3]{Lan} $\oint:  \bbS_{K,\Sigma} \ra \bbS_K^{min}$.    For each standard maximal rational parabolic $P \subset G$, we let the $P$-stratum $\partial^P\bbS_{K,\Sigma}$ of the toroidal boundary of $\bbS_{K,\Sigma}$ be the closure of $\oint^{-1}(\partial^P\bbS_K^{min})$.  Let $R \subset G$ be a standard rational parabolic, and write $R = \cap Q_j$ for a (unique) set of
standard maximal parabolics $Q_j$.  We define the $R$-stratum: 
\begin{equation}\label{Rstratum}
\partial^R\bbS_{K,\Sigma} = \cap_j \partial^{Q_j}\bbS_{K,\Sigma}.
\end{equation}
Under our running hypotheses on $\Sigma$, the intersections are componentwise transversal and each $\partial^R\bbS_{K,\Sigma}$
is a union of smooth subvarieties of codimension equal to the parabolic rank $r(R)$ of $R$, intersecting transversally.  

Suppose $R$ has parabolic rank $r$.  We define the nerve $\mathfrak{N}_{\Sigma}(R)$ of the closed covering
of $\partial^R\bbS_{K,\Sigma}$ by  irreducible components of codimension $r$, 
as in \cite[\S 3.1]{HZ01} (where the $R$-stratum was denoted $Z_\Sigma(R)$).   We recall the relation between the homotopy type of
$\mathfrak{N}_{\Sigma}(R)$ and  Borel-Serre compactifications in  \S \ref{homotopynerve}.

%\subsubsection{Topological properties of $\Sigma_F$}\label{nonmax}  This section collects topological facts about  The reader is advised to skip this section
%until these facts are needed.

Let $R$ be a standard rational proper parabolic subgroup of $G$, and define the maximal parabolic $P(R)$ as in \S \ref{parab}; we suppose $P(R) = P = P_F$, as above, and define $F(R) = F$.   We can write $R = \cap_{j = 1}^{r(R)} Q_j$ uniquely as the intersection of  standard rational maximal parabolics,
with $Q_j < Q_{j'}$, in the total order defined above, if and only if $j < j'$.  Thus $Q_1 = P(R)$; $P(R)$ is the smallest maximal parabolic containing $R$, in the order by size of $G_h$, but it has the largest $G_\ell$.   The cone $C_F$ is homogeneous under the action of 
$G_{\ell,P(R)}(\RR) = GL(m(R),\CC)^d$, and can be identified with the symmetric space attached to the Lie group $GL(m(R),\CC)^d$.  The partial compactification
$\bar{C}_F$ is a subset of the topological closure of $C_F$ in $U_F(\RR)$ that is set-theoretically the union of cones of the form $C_{F'}$, where $F'$ runs over boundary components of $X$ such that the stabilizer $P_{F'}$ contains $G_{h,R}$; alternatively, such that $F$ is a boundary component of $F'$.  Here it is necessary to allow $P_{F'}$ to be a non-standard maximal parabolic; the construction of the toroidal compactification involves canonical rational embeddings $U_{F'} \subset U_F$ for all such $F'$.
Alternatively, we can take the union of $C_{F'}$ for which $P_{F'} > P(R)$, together with their translates under $G_{\ell,P(R)}(\QQ) \simeq GL(m(R),F)$ (here
$F$ designates the CM field!).    

Recall that $\Sigma_F$ is a fan in $\bar{C}_F$, and in particular can be written as a disjoint union 
\begin{equation}\label{Ffans} \Sigma_F = \coprod_{F'} \Sigma_F(F')
\end{equation}
where $\Sigma_F(F')$ is the set of polyhedral cones in $\Sigma_F$ contained in $C_{F'}$.   To fix ideas we assume for the moment that 
$\Gamma_{\ell,P} \subset G_{\ell,P}(\QQ)$ is a congruence subgroup and that $\Sigma_F$ is invariant under $\Gamma_{\ell,P}$; this is one of the properties
of $\Sigma$ that is used to construct classical (non-adelic) toroidal compactifications in \cite{AMRT}.  We assume $\Gamma_{\ell,P}$ is neat.  Then the locally symmetric space $X(\Gamma_{\ell,P}) := \Gamma_{\ell,P} \backslash C_F$ has a Borel-Serre compactification $\bar{X}(\Gamma_{\ell,P})$ and the inclusion
$X(\Gamma_{\ell,P}) \subset \bar{X}(\Gamma_{\ell,P})$ is a homotopy equivalence.    

The relation between toroidal
and Borel-Serre compactifications  will be needed for the computations in \S \ref{topologicalpart}; see \S \ref{homotopynerve}.

%{\color{red}

%Introduce, either here or in \S \ref{topologicalpart}, a general statement about the cohomology of Borel-Serre strata that replaces Kostant's
%formula.
%}
%%%%%%%%%%%%%%%%%%%%%%%%%%%%%%%%%%%%%%%%%%%%%%%%%%%%%%%%%%%%%%%%%%%%%%%%%%%%%%%%%%%%%%%%
%%%%%%%%%%%%%%%%%%%%%%%%%%%%%%%%%%%%%%%%%%%%%%%%%%%%%%%%%%%%%%%%%%%%%%%%%%%%%%%%%%%%%%%%

\subsection{Level subgroups}

We will be working in the following situation.  Start with $K$ as above, and let $S = S(K)$ be its set of ramified primes, as in \S \ref{avbnotation}.  
We recall that $S$ contains a prime $\frr$ with a divisor $\frr_0$ and that $K_{\frr_0}$ is chosen to guarantee
that $K$ is neat.  Choose a finite set $Q$ of primes of $F^+$ that split in $F$,  with the property that, if $v \in Q$ divides the rational prime $q$,
then $q \notin S$ and $q \equiv 1 \pmod{p}$.    Let $S(Q)$ be the set of rational primes divided by primes in $Q$; we assume that each $q \in S(Q)$ is divided by a unique $v \in Q$.  
For each $v \in Q$ let $\fm_v \in \CO_v$ be the maximal ideal, $k(v) = \CO_v/\fm_v$ the residue field, and let 
$$K_{0,v} = \big\{k \in GL(n,\CO_v) ~|~ k \equiv \begin{pmatrix} a(k) & b \\ 0 & d \end{pmatrix} \pmod{\fm_v}\big\}$$
with $a(k) \in k(v)^{\times}$, $b$ a $1 \times (n-1)$ row matrix over $k(v)$, and $d \in GL(n-1,k(v))$.  Let $K_{1,v} \subset K_{0,v}$ denote the kernel of the map $k \mapsto a(k)$, and let
$K_{\Delta,v}$ be the smallest subgroup of $K_{0,v}$ containing $K_{1,v}$ such that $\Delta_v = K_{0,v}/K^{\Delta}_v$ is a $p$-group.  Let $\Delta_Q = \prod_{v \in Q} \Delta_v$.
For $q \in S(Q)$ let 
$K_{i,q} = K_{i,v} \times [\prod_{v' \in D_+(q)\setminus v} K_{v'}] \times K_{-,q}$, for $i = 0, 1, \Delta$, and let
$$K_{i,Q} = \prod_{q \notin S(Q)} K_q \times \prod_{q \in S(Q)} K_{i,q}.$$

Let $w \notin S$ be a rational prime, fix $t \in G_w$, and let $K_w(t) = K_w \cap t K_w t^{-1}$.
Suppose $w \notin S \coprod S(Q)$, and $t \in G_w$.  Define the {\it modification} $K_{i,Q}(t)$ of $K_{i,Q}$, with $i = 0, 1, \Delta$, to be the product
$$K_{i,Q}(t) = \prod_{q \notin S(Q)\cup w} K_q \times \prod_{q \in S(Q)} K_{i,q} \times K_w(t).$$
If $Q$ is empty, one writes $K(t) = K_{0,Q}(t) = K_{1,Q}(t)$.  

The following theorem is proved in  \S A.6 of the first author's thesis \cite{At}.

\begin{thm}\label{etale}  Suppose $K=\prod_q K_q\subseteq G(\af)$ is a neat subgroup, such that $K_q$ are maximal hyperspecial for almost all $q$. Suppose further that $K\supset K_p$, a fixed hyperspecial maximal $K_p\subseteq G(\QQ_p).$ Fix a set $Q$ of primes as above.  Let $\CB$ be the set of modifications $K_{i,Q}(t)$ of $K_{i,Q}$, with $i = 0, 1, \Delta$.  
Then for every $K_{i,Q}(t)$ there is a toroidal datum
$\Sigma(K_{i,Q}(t))$ such that
\begin{itemize}  
\item[(1)] The toroidal compactification $S_{K_{i,Q}(t)}(G,X) \hookrightarrow   S_{K_{i,Q}(t)}(G,X)_{\Sigma_{i,Q}(t)}$ is smooth and projective;
\item[(2)]  If $K' \subset K$ is an inclusion in $\CB$ then the natural covering map 
$$S_{K'}(G,X) \ra S_{K}(G,X)$$ extends to a map 
$$S_{K'}(G,X)_{\Sigma(K')} \ra S_{K}(G,X)_{\Sigma(K)}$$ of toroidal compactifications. 
\item[(3)] The map 
$$S_{K_{\Delta,Q}}(G,X)_{\Sigma_{\Delta,Q}} \ra S_{K_{0,Q}}(G,X)_{\Sigma_{0,Q}}$$ 
is an \'etale covering with group 
$\Delta_Q$.
\item[(4)]  Conditions (1)-(3) hold when the $S_{K_{i,Q}}(G,X)$ are replaced by their integral models $\bbS_{K_{i,Q}}$.
\end{itemize}
\end{thm}
The first two assertions follow from the classical theory of toroidal compactifications and are used without comment in defining Hecke correspondences in \S \ref{Heck-op}. The third one requires more delicate analysis of the analytic construction in \cite{AMRT}.  Atanasov uses the observation that the intersection of $U_F(\QQ_q)$ with any compact open subgroup is contained in a pro-$q$ group for each $F$ and each unramified prime $q \neq p$.   Since $\Delta_Q$ is a $p$-group, it then suffices to show that $K_{\Delta,Q}$ and $K_{0,Q}$ restrict to the same subgroup of $U_F(\QQ)$, and this is a simple calculation.  The last claim follows by purity of the branch locus applied to each irreducible boundary component, which reduces it to the characteristic zero statement in (3).

\subsection{Automorphic vector bundles over toroidal compactifications}

In the next few sections we will develop the theory of coherent cohomology and Hecke algebras without controlling
ramification at the places in $S \setminus S(G)$, in the notation of \S \ref{avbnotation}, in order to avoid unnecessary
complications.   Sections \ref{ramtypes} and \ref{typedata} will indicate the adjustments that are needed in order to accommodate ramification.

Fix an automorphic vector bundle $\CF = \CF_{\tau}$ on $S_{K_p}(G,X)$, and identify $\CF$ with the family of $\CF_K$, as in \ref{avbnotation}.  
For each open compact $K \supset K_p$, let $\Sigma_K$ be a cone decomposition defining a toroidal compactification $S_K(G,X) \hookrightarrow S_K(G,X)_{\Sigma_K}$.   Provided $K$ is neat, each $\CF_K$ admits a {\it canonical extension} $\CF_K^{can}$ over $S_{K}(G,X)_{\Sigma_K}$, constructed in \cite{H89} following Mumford.   The canonical extensions are functorial with respect to pullback and direct image:  if $K' \subset K$ is an open subgroup, and if $\Sigma_{K'}$ is a cone decomposition for $S_{K'}(G,X)$ compatible with $\Sigma_K$, so that the natural finite covering map 
$$\pi_{K,K'}:  S_{K'}(G,X) \ra S_K(G,X)$$
extends to a morphism of toroidal compactifications
$$\pi_{K,K'}^{tor}:  S_{K'}(G,X)_{\Sigma_{K'}} \ra S_K(G,X)_{\Sigma_{K'}}$$
then there are canonical isomorphisms
$$\pi_{K,K'}^{tor,*}(\CF_K^{can}) \isoarrow \CF_{K'}^{can}~~\text{and}~~(\pi_{K,K'}^{tor})_*(\CF_{K'}^{can}) \isoarrow \CF_{K}^{can},$$
and moreover
$$ R^i(\pi_{K,K'}^{tor})_*(\CF_{K'}^{can}) = 0, ~~ i > 0.$$

Note that this applies in particular when $K' = K$ but $\Sigma_{K'}$ is a refinement of $\Sigma_K$.  

Similarly, suppose the complement of $S_K(G,X)$ in  $S_K(G,X)_{\Sigma_K}$ is a divisor $D_{\Sigma_K}$ with normal crossings.  Let  $\mathcal{I}_{\Sigma_K}$ be the ideal sheaf of $D_{\Sigma_K}$ -- it is a line bundle on $S_K(G,X)_{\Sigma_K}$ -- and let $\CF_K^{sub} = \CF_K^{can}\otimes\mathcal{I}_{\Sigma_K}$ be the {\it subcanonical extension} of $\CF_K$.   Then, under the hypotheses above, there is a canonical isomorphism
$$(\pi_{K,K'}^{tor})_*(\CF_{K'}^{sub}) \isoarrow \CF_{K}^{sub},\quad R^i(\pi_{K,K'}^{tor})_*(\CF_{K'}^{sub}) = 0, ~~ i > 0.$$
Moreover, 
\begin{equation}\label{pullbacketale} \text{ if $\pi_{K,K'}^{tor}$ is \'etale, then one also has } \pi_{K,K'}^{tor,*}(\CF_K^{sub}) \isoarrow \CF_{K'}^{sub}. \end{equation}

Except in a few low-dimensional
cases, the action of $G(\afp)$ only extends to the family $S_{K,\Sigma}(G,X)$ if the $\Sigma$ are also allowed
to vary along with $K$, and in particular the algebra $\CH_K$
of Hecke operators of level $K$ do not act geometrically on $\CF_K^{can}$.
However, it is explained in \cite{H90} that $\CH_K$ does act canonically on the finite-dimensional vector
space $H^i(S_{K,\Sigma},\CF_K^{can})$, for each $i$.  We reformulate the result of \cite{H90}, and the integral version proved
in \cite{Lan3}, in a version better adapted to localization of Hecke algebras at maximal ideals.
Let $E(\CF)$ denote an extension of the reflex field $E(G,X)$ over which the automorphic vector bundle $\CF$ has a rational model.  Most of the following Proposition is a consequence of the properties of canonical and subcanonical extensions recalled above.

\begin{prop}\label{db}
  \begin{itemize}
\item[(i)]  Let $(G,X)$, $K$, and $\Sigma_K$ be as above.   Let $\Sigma'_K$ be a refinement of $\Sigma_K$.  Then the morphism 
$$\pi_{\Sigma,\Sigma'}:  S_{K}(G,X)_{\Sigma'_{K}} \ra S_K(G,X)_{\Sigma_{K}}$$
induces (by pushforward) a canonical quasiisomorphism
$$R\pi_{\Sigma,\Sigma',*}:  R\Gamma(S_K(G,X)_{\Sigma'_{K}}, \CF_K^{can}) \simeq R\Gamma(S_K(G,X)_{\Sigma_{K}}, \CF_K^{can})$$
in the bounded derived category $D^b(E(\CF))$ of complexes of $E(\CF)$-vector spaces.  (Here the superscript $^{can}$ is used to designate the canonical extension for any toroidal compactification.)
Similarly, $\pi_{\Sigma,\Sigma'}$ induces a canonical quasiisomorphism
$$R\pi_{\Sigma,\Sigma',*}:  R\Gamma(S_K(G,X)_{\Sigma'_{K}}, \CF_K^{sub}) \simeq R\Gamma(S_K(G,X)_{\Sigma_{K}}, \CF_K^{sub})$$

In particular, the objects 
$$R\Gamma(S_K(G,X)_{\Sigma_{K}}, \CF_K^{can}), ~~R\Gamma(S_K(G,X)_{\Sigma_{K}}, \CF_K^{sub})$$
 of $D^b(E(\CF))$  are well-defined and independent of the choice of $\Sigma_K$; we denote them
$R\Gamma^{can}(S_K(G,X),\CF_K)$ and $R\Gamma^{sub}(S_K(G,X),\CF_K)$, respectively.

\item[(ii)]  Let $t \in G(\af)$.  For any $K$ as above, and for any $\CF$, the natural isomorphism
$$*t^{-1}:  S_K(G,X) \ra S_{tKt^{-1}}(G,X)$$
defines an isomorphism in $D^b(E(\CF))$ by pullback:
$$*t^{-1}: R\Gamma^{can}(S_{tKt^{-1}}(G,X),\CF_{t^{-1}Kt})  \simeq R\Gamma^{can}(S_K(G,X),\CF_K).$$
There are analogous isomorphisms when $^{can}$ is replaced by $^{sub}$.

\item[(iii)]  Let $K' \subset K$ be an open subgroup.  The pullback and direct image functors define canonical morphisms in $D^b(E(\CF))$:
$$\pi_{K,K'}^*  R\Gamma^{can}(S_K(G,X),\CF_K) \ra R\Gamma^{can}(S_{K'}(G,X),\CF_{K'});$$
$$\pi_{K,K'}^*  R\Gamma^{sub}(S_K(G,X),\CF_K) \ra R\Gamma^{sub}(S_{K'}(G,X),\CF_{K'});$$
$$\pi_{K,K',*}  R\Gamma^{can}(S_{K'}(G,X),\CF_{K'}) \ra R\Gamma^{can}(S_{K}(G,X),\CF_{K})$$
$$\pi_{K,K',*}  R\Gamma^{sub}(S_{K'}(G,X),\CF_{K'}) \ra R\Gamma^{sub}(S_{K}(G,X),\CF_{K})$$

\item[(iv)]  Let $\CO \subset E(G,X)$ denote the localization at a prime above $p$ of the integers of $E(G,X)$.  With $K$ and $K'$ as above, both assumed unramified at $p$, suppose $S_K(G,X)$, $\CF_K$, $S_{K'}(G,X)$, $\CF_{K'}$, and their toroidal compactifications all admit compatible integral models over $\CO$, in the sense of Proposition \ref{toroidalintegral}.  Denote these models $\bbS_K$, 
$\bbS_{K,\Sigma_K}$, $\CF_K$ (no change), etc.  Then the object  
$R\Gamma(\bbS_{K,\Sigma_{K}}, \CF_K^{can})$ of $D^b(\CO)$ is a well-defined perfect complex of $\CO$-modules and independent of the choice of $\Sigma_K$; we denote it 
$RH(\CF_K)$, or when necessary, $RH(\CF_K^{can})$.  The conclusions of (iii) hold for the corresponding elements of $D^b(\CO)$.   If $t \in G(\afp)$, then the conclusions of (ii) holds for the corresponding elements of $D^b(\CO)$.  

The analogous statements hold when $^{can}$ is replaced by $^{sub}$. 

\end{itemize}

\end{prop}

\begin{proof}  Items (i)-(iii) are proved in \cite{H90}; item (iv) is due to Lan and is proved in Proposition 1.4.3 of \cite{Lan3}.
\end{proof}

\section{Hecke algebras}\label{Heckealgebras}   
\subsection{Adelic representations and fixed vectors}

Let $K = K_p \times K^p$, where $K_p$ is a fixed hyperspecial maximal compact subgroup of $G(\QQ_p)$ as before and $K^p$ is an open compact subgroup of $G(\afp)$.   We define
$$RH(\CF_{K_p}) = RH(\CF_{K_p}^{can}) = \varinjlim_{K^p} RH(\CF_{K_p\times K^p}^{can})$$
where the limit is taken as $K^p$ shrinks to the identity.    This is a well-defined object of the bounded derived category of $\CO$-modules, and since the colimit is exact on this category we have
\begin{equation}\label{colcoh}
H^i(RH(\CF_{K_p}^{can})) \isoarrow \varinjlim_{K^p} H^i(RH(\CF_{K_p\times K^p}^{can})
\end{equation}
for all $i \in \ZZ$.

We define $RH(\CF_{K_p}^{sub})$ similarly and note that the analogue of \eqref{colcoh} is valid with $can$ replaced by $sub$.  

Proposition \ref{db} (ii) and (iv) have the following Corollary:
\begin{cor}\label{db1}  The objects $RH(\CF_{K_p}^{can})$, $RH(\CF_{K_p}^{sub})$ have canonical compatible actions of $G(\afp)$.  Moreover, the actions are smooth in the sense that any element of $H^i(RH(\CF_{K_p}^{?}))$, where $? = can, sub$ and $i \in \ZZ$, is fixed by an open compact subgroup of $G(\af)$.
\end{cor}

The last sentence follows from \eqref{colcoh}.  However, it is important to note that, for a given open $K \supset K_p$, it is not obviously the case that the canonical map
\begin{equation}\label{fixedK} H^i(\bbS_{K,\Sigma_{K}}, \CF_K^{?}) \rightarrow H^i(RH(\CF_{K_p}^{?}))^K
\end{equation}
is an isomorphism, or even necessarily injective.   Instead there is a Hochschild-Serre spectral sequence
\begin{equation}\label{HoSe}
E_2^{r,s} = H^r(K,R^sH(\CF_{K_p}^?)) \Rightarrow H^{r+s}(\bbS_{K,\Sigma_{K}}, \CF_K^{?}),
\end{equation}
again with $? = can, sub$.  However, we have the following result:

\begin{prop}\label{BCGP}\cite{BCGP1} With notation as above,  $H^r(K,R^sH(\CF_{K_p}^?)) = 0$ for all $r > 0$, $? = can, sub$.  Thus the canonical map 
$$RH(\CF_{K}^{?}) \rightarrow RHom_{K^p}(1,RH(\CF_{K_p}^{?})),$$
where $1$ denotes $\CO$ with the trivial action of $K^o_\Psi$, is a quasi-isomorphism for $? = can, sub$.
\end{prop}

The Proposition is proved in \cite{BCGP1} by verifying that the stabilizers in $K^p$ of points in the toroidal boundary have order prime to $p$.  

Let $S = S(K)$ be the set of primes where $K$ does not contain a hyperspecial maximal compact subgroup, plus the prime $p$, and let $G(\afS)$ be the set of ad\`eles of $G$ with entry $1$ at primes in $S$; let $K_S = K\cap G(\afS)$.  Let $\CO$ be as above and define the Hecke algebra $\CH_{K,\CO}$ to be the $\CO$-algebra of compactly-supported $K_S$-biinvariant functions on $G(\afS)$, with multiplication given by convolution.  Since $K_S$ contains a hyperspecial maximal compact subgroup at each prime, $\CH_{K,\CO}$ is a commutative algebra.  Corollary \ref{db1} and Proposition \ref{BCGP} thus imply that $\CH_{K,\CO}$ acts on 
$RH(\CF_{K}^{?}))$ for $? = can, sub$.  We le $T_{S,K} = T_{S,K,\CF}$ denote the $\CO$-subalgebra of $\Endo(RH(\CF_K))$ generated by this action.

Since \cite{BCGP1} is not yet available, we give an alternative construction of $T_{S,K}$ in the next section.

\subsection{Hecke operators}\label{Heck-op}

Notation is as in Proposition \ref{db}.   Fix $K$ hyperspecial at $p$ and let $S = S(K)$ be as above.  %Let $\Psi$ be a finite set of primes $q$ not in $S \cup \{p\}$ and write $K^p = K^{p,\Psi} \times K_\Psi$, where $K_\Psi \subset G_\Psi = \prod_{q \in \Psi} G(\QQ_q)$ and $K^{p,\Psi}$ is an open compact subgroup of the ad\`eles of $G$ away from $\Psi \cup \{p\}$.  Let $K_{\Psi}^o \subset G_\Psi$ be a (product over $\Psi$ of ) hyperspecial maximal compact subgroups.  For $? = can, sub$, we define a variant of $RH(\CF_{K_p}^{?})$:
%$$RH(\CF_{K_p,\Psi}^{?}) = \varinjlim_{K_\Psi} RH(\CF_{K_p\times K^{p,\Psi} \times K_\Psi}^{?})$$
%where $K_\Psi$ shrinks to the identity.  Then $G_\Psi$ acts on $RH(\CF_{K_p,\Psi}^{?}$ and therefore the $\CO$-subalgebra $T_\Psi$ of the group algebra $\CO[G_\Psi]$ of $K^o_\Psi$-biinvariant functions acts on
%\begin{equation}\label{rhom}
%RH(\CF_{K_p,\Psi}^{?})^{K^o_\Psi} := RHom_{K^o_\Psi}(1,RH(\CF_{K_p,\Psi}^{?})),
%\end{equation}
%where $1$ denotes $\CO$ with the trivial action of $K^o_\Psi$.  
%{\color{red}.   with the property that all primes of $F^+$
%dividing $q$ split in $F/F^+$.  }
For $t \in G(\afS)$, a Hecke operator 
$T(t) \in \Endo(RH(\CF_{K}^{?}))$, denoted $T_t^\Sigma$, is defined in \cite[\S 8.1.5]{GK}, for $? = can, sub$.  An alternative definition is given in \cite[\S 4.2]{BP}, where it is also proved that the subalgebra $T_{S,K}$ of $\Endo(RH(\CF_{K}^{?}))$ generated by the $T(t)$ is commutative.
%let $K(t)$ be as above, and consider the composition of morphisms in $D^b(\CO)$:
%\begin{multline}\label{hecket}
%R\Gamma^{can}(\bbS_K,\CF_K) \overset{\pi_{K,K(t)}^*}{\longrightarrow} R\Gamma^{can}(\bbS_{K(t)},\CF_{K(t)}) \\ \overset{\pi_{t^{-1}Kt,K(t),*}}{\longrightarrow} R\Gamma^{can}(\bbS_{t^{-1}Kt},\CF_{t^{-1}Kt})
%\overset{*t^{-1}}{\longrightarrow} R\Gamma^{can}(\bbS_K,\CF_K).
%\end{multline}
%Denote the composition $T(t)$; it is an element of $\Endo(R\Gamma^{can}(\bbS,\CF_K))$.  Under the hypotheses of (iv) of Proposition \ref{db}, $T(t)$ defines an element of $\Endo(RH(\CF_K))$.   Let $S = S(K)$ as above, and let $T_{S,K} = T_{S,K,\CF}$ denote the $\CO$-subalgebra of $\Endo(RH(\CF_K))$ generated by the $T(t)$ for $t \in G_q$, where $q$ runs through the set of primes  not in $S \cup \{p\}$ 
%It is a commutative algebra, because the Hecke operators at places $q$ at which $K_q$ is hyperspecial commute with one another.   
The algebra $T_{S,K}$ naturally maps to $\Endo(\oplus_i H^i(RH(\CF_K)))$ but the map to its image is not necessarily an isomorphism.  In particular, there is no reason to assume $T_{S,K}$ to be reduced.  %{\color{red}[May prefer to take only $t \in G_q$ when all primes above $q$ are split in $F/F^+$, since this always seems to be the hypothesis in the literature.]}

In the setting of (iv) of Proposition \ref{db}, the natural inclusion $\CF_K^{sub} \ra \CF_K^{can}$ (the cone decomposition is implicit and omitted from the notation) determines a morphism in $D^b(\CO)$
\begin{equation}\label{cone} R\Gamma^{sub}(\bbS_K,\CF_K) \ra R\Gamma^{can}(\bbS_K,\CF_K) \end{equation}
Let $R\Gamma^{\partial}(\bbS_K,\CF_K)$ denote the cone on the morphism \eqref{cone}.  

The following theorem summarizes the state of the art:

\begin{thm}\label{Galoisreps}  (i)  Let $\kappa$ be an algebraically closed field and let $\nu:  T_{S,K,\CF} \ra \kappa$ be a continuous homomorphism.  Then there is a semisimple $n$-dimensional representation
$\rho_{\nu}:  \Gamma_{F} \ra GL(n,\kappa)$ that is characterized, up to equivalence, by the following property:
%\begin{itemize}
%\item[(a)]  
If $v$ is a prime of $F$ not in $S  \cup \{p\}$, then $\rho_{\nu}$ is unramified at $v$.  Let $\Gamma_v \subset \Gamma_F$ be a decomposition group at $v$; then the semisimplification of the restriction $\rho_{\nu,v}$ of $\rho_{\nu}$ to $\Gamma_v$ corresponds to the restriction to $\nu$ to the image of the Hecke operators at $v$ by the unramified Langlands correspondence, in the following sense:  there is an equality of polynomials in $\kappa[X]$
\begin{equation}\label{charpoly}
\det(1 - \rho_{\nu}(\Frob_v)X) =   1 + \sum_{i = 1}^n (-1)^i \nu(q^{\frac {(n+1)i}{2}}T_{i,v} )X^i
\end{equation}
where the $T_{i,v}$ are the standard Hecke operators at $v$, normalized as in \cite{Sch}. %?? 

%\item[(b)]  Suppose $\kappa$ is of characteristic zero and $\nu$ is the homomorphism defined by the $K$-fixed vectors of a cuspidal cohomological automorphic representation $\pi$ of $GL(n)$, or of a cuspidal  and $v \in S$, the restriction of $\rho_{\nu}$ to $\Gamma_v$ 

%\end{itemize}

(ii)  Suppose $\kappa$ is of characteristic zero and $\nu$ is the homomorphism attached to a cuspidal automorphic representation of $G$.  Then there is a homomorphism
$$r_{\nu}:  \Gamma_{F^+} \ra \CG_n(\kappa)$$
where $\CG_n$ is the disconnected algebraic group defined in (\cite{CHT}, \S 2.1) that corresponds to the pair $(\rho_{\nu},\beta)$, for some character $\beta$ of $\Gamma_F$, under the correspondence of Lemma 2.1.1 of \cite{CHT}.  In particular,
the identity component $\CG_n^0$ is isomorphic to $GL(n) \times GL(1)$, and the restriction of $r_{\nu}$ to $\Gamma_F \subset \Gamma_{F^+}$ is given by $(\rho_{\nu},\beta)$.

(iii)  Write $R\Gamma^{\partial}(\CF_K) := R\Gamma^{\partial}(\bbS_K(G,X),\CF_K)$.  Suppose $\nu$ occurs in the representation on some subquotient of $R\Gamma^{\partial}(\CF_K)$; we say that $\nu$ occurs in the support of $R\Gamma^{\partial}$ (for some $\CF_K$).  Then $\rho_{\nu}$ is reducible.
\end{thm}

Part (ii) is the familiar association of Galois representations to polarized cohomological cuspidal automorphic representations of $GL(n)$.  Part (i) includes the results of \cite{HLTT} and the extension to torsion cohomology by Scholze, Boxer, Pilloni-Stroh, and Goldring-Koskivirta \cite{Sch,Box,PS,GK}.  

  It remains to prove Part (iii), which follows from part (i) and an analysis of the non-cuspidal coherent cohomology of $\bbS_K(G,X)$, based primarily on the considerations of \cite{HZ94}, \cite{LS2}, and \cite{Lan2}.  The details are postponed until section \ref{proofofreps}; the proof will be particularly long, because it requires a review of the structure of the toroidal boundary.

Let $\nu$ be as in Theorem \ref{Galoisreps}, with $k$ a field of characteristic $p$.   Let $\fm = \ker{\nu}$; we write $\fm = \fm_{\nu}$ or $\fm = \fm_{\rho}$ with $\rho = \rho_{\nu}$.  Let 
$\TT_{\nu}$ be the localization of $T_{S,K}$ at $\fm_{\nu}$, and let 
$$RH(\CF_K)_{\nu} = RH(\CF_K)\otimes_{T_{S,K}}\TT_{\nu}$$ denote the localization of $RH(\CF_K)$ at $\fm_{\nu}$.   
Recall that $(\tau,W_\tau)$ is the representation of the parabolic $P_x$ corresponding to a point $x \in \check{X}$.   Let $\mu_\tau \in X^+_x(T)$ denote the highest weight of $\tau$, and assume 
\begin{equation}\label{mumu}
\mu_\tau = w\cdot \mu, w \in W^x, \mu \in X^+(T).
\end{equation}

\begin{prop}\label{torsionfree}  We assume $\CF_K$ is the automorphic vector bundle attached to a character $w\cdot\mu$,
as in \eqref{mumu}.
Suppose (a) $\rho_{\nu}$ is irreducible; (b) the highest weight $\mu_\CF$
of $\CF$ satisfies the Lan-Suh vanishing conditions:

\begin{itemize}
\item[(i)]  The character $\mu$  of \eqref{mumu} is $p$-small (Definition \ref{psmall});
\item[(ii)] $\mu$ is {\it sufficiently regular} in the sense of \cite[7.18]{LS1}
\item[(iii)] $\mu$ satisfies the inequality \cite[7.22]{LS1}.
\end{itemize}

Then  the cohomology of $RH(\CF_K)_{\nu}$ is concentrated in a single degree $i_0 = i_0(\CF)$ and is torsion-free over $\CO$.   
\end{prop}

\begin{proof}  Write $\bbS = \bbS_K(G,X)$.  For $q \geq 0$, define the  {\it interior cohomology} of $\CF_K$
to be
$$H^q_!(R\Gamma(\bbS,\CF_K)) :=  \textrm{Im}[H^q(\bbS,\CF_K^{sub}) \ra H^q(\bbS,\CF_K^{can})].$$
It follows from hypothesis (a) and (iii) of Theorem \ref{Galoisreps} that the localization at $\fm_{\nu}$ of $R\Gamma^{\partial}(\bbS,\CF_K)$ is acyclic.  
Thus for any $q$, 
$$H^q(\bbS,\CF_K^{can})_{\nu} = H^q_!(R\Gamma(\bbS,\CF_K))_{\nu};$$ 
in other words, the localization at $\fm_{\nu}$ of the cohomology of the canonical extension coincides with the localization at $\fm_{\nu}$ of the  interior cohomology.  In view of hypothesis (b), it follows from  Corollary 7.25 and Theorem 8.2 (3) of \cite{LS2} that $H^q_!(R\Gamma(\bbS,\CF_K))_\nu$ is concentrated in a single degree $i_0(\CF)$ and is torsion-free over $\CO$.  The Proposition is an immediate consequence.
\end{proof}

The Lan-Suh vanishing conditions are written out more explicitly in \cite[\S 6.10]{H13}.
We write $H^q_!(\bbS,\CF_K)$ instead of $H^q_!(R\Gamma(\bbS,\CF_K))$.  

\subsection{Ramification and types}\label{ramtypes}  The modules $H^q_!(R\Gamma(\bbS,\CF_K))$  of coherent cohomology
give rise to Galois representations that may be ramified  at places $v$ where $K_v$ is not hyperspecial maximal.  
In order to restrict attention to minimal deformation problems, we introduce spaces of cohomology with coefficients in {\it semisimple Bushnell-Kutzko types} \cite{BK99}.  
This was explained in the unpublished manuscript \cite{HT}, in the case of supercuspidal types, but seems have never been used subsequently, although the definition we give below of automorphic forms with coefficients in Bushnell-Kutzko types  is essentially the one that appeared in \cite{HT} as well as \cite{V,K01}; for automorphic forms on
$GL(2)$ it has been standard for some time.  A version based on finite-dimensional representations of 
the multiplicative groups of local division algebras is developed in \cite[p. 97]{CHT}; this is essentially equivalent for purposes of proving modularity but it cannot be applied
to prove that the coherent cohomology is free over the localized Hecke algebra, which is our main theorem.  

Let $S = S(K)$ be as above, and assume the local groups $K_v$ and $K_{-,q}$ are as in \S \ref{avbnotation}.   Let $(K_v^+,K_v,\Yu_v,\Lambda_{\Yu_v})$ be type data as in 
\S \ref{typedata}.
\begin{defn}\label{semisimpletype}  The pair $(K_v^+,\Yu_v)$ is a semisimple type if it 
determines a component of the Bernstein center of $GL(n,F^+_v)$ in the following sense:
\begin{itemize}
\item[(a)]  For any irreducible representation $\pi$ of $GL(n,F^+_v)$, 
$$\dim \Hom_{K_v^+}(\Yu_v,\pi) \leq 1;$$
\item[(b)]  If $\pi$ is an irreducible representation of $GL(n,F^+_v)$
then 
$$\Hom_{K_v^+}(\Yu_v,\pi) \neq 0$$
 if and only if there is a parabolic subgroup $P \subset GL(n,F^+_v)$ with Levi quotient $L$, 
a (fixed) supercuspidal representation $\sigma$ of $L$, and an unramified character $\alpha$ of $L$, such that $\pi$ is an irreducible constituent
of the  induced representation $Ind_P^{GL(n,F^+_v)} [\sigma \otimes \alpha]$.
\end{itemize}
\end{defn}

Let $\Yu = (\Yu_v)_{v \in S}$, where for all $v \in S$ (we are now using $S$ to denote a set of primes of $F^+$, as in \S \ref{avbnotation}) we assume 
\begin{hyp}\label{Yuhyp}

(1) Either  $\Yu_v$ is 
the trivial representation of the special maximal compact subgroup $K_v$, in which case we set $K_v^+ = K_v$, or the pair $(K_v^+,\Yu_v)$ is a semisimple type for the irreducible representation $\pi$
of $GL(n,F^+_v)$.    Let $S^+ \subset S$ be the subset for which the second condition holds, and let $K^+_S = \prod_{v \in S} K_v^+$.

(2)  For all $v \in S^+$, $p$ is a banal prime for $GL(n,F^+_v)$:  that is, $p$ does not divide the pro-order of any compact open subgroup
of $GL(n,F^+_v)$.
\end{hyp}

We let $\Lambda_{\Yu}$ denote the representation $\otimes_v \Lambda_{\Yu_v}$ of $\prod_{v \in S} K_v^+$.
Hypothesis \ref{Yuhyp} (2) implies that the module $\Lambda_{\Yu}$ is the unique $\prod_{v \in S} K_v^+$-invariant lattice in $\Lambda_{\Yu_v}\otimes \QQ$,
up to scaling.   With $\CF_K$ an automorphic vector bundle as above, we let 
$$\CF_{K,\Yu} = \CF_K\otimes_\CO \Lambda_{\Yu}$$ and for any cohomology module $\mathbf{H}$ with coefficients in $\CF_{K,\Yu}$ we define
\begin{equation}\label{cohoYu}
\mathbf{H}_{\Yu} \coloneqq \Hom_{K_S^+}(\Lambda_{\Yu}, \mathbf{H}).
\end{equation}
Thus we can define $H^q_!(R\Gamma(\bbS,\CF_{K,\Yu}))$, $H^q(\bbS,\CF_{K,\Yu}^{can})_{\nu,\Yu}$, etc.
\smallskip

\noindent{\bf Theorem \ref{Galoisreps}}  {\it  (iv) The assertions of Theorem \ref{Galoisreps} (iii)  hold when $\CF_K$ is replaced by $\CF_{K,\Yu}$.}

%We will omit the subscript $\Yu$ and write  $T_{S,K,\CF}$ for the subalgebra of $\Endo(RH(\CF_{K,\Yu}))$ defined as in \S \ref{Heckealgebras}.  

\subsection{Minimality condition}\label{minimality}
For applications to the Taylor-Wiles method, we add the following hypotheses on the Galois representations attached to the type data.  For each $v$, we choose
a $\pi$ in the inertial equivalence class corresponding to the quadruple $(K_v^+,K_v,\Yu_v,\Lambda_{\Yu_v})$ -- in other words,
$\Hom_{K_v^+}(\Yu_v,\pi) \neq 0$.  Let $\cL_v(\pi)$ denote the $n$-dimensional representation of the Weil-Deligne group $W_v$ of $F^+_v$ attached
to $\pi$ by the local Langlands correspondence, with coefficients in the $p$-adic integer ring $\CO$, and let $\overline{\cL}_v(\pi)$ denote its reduction
modulo  the maximal ideal $\fm_\CO$.  We attach to $\cL_v(\pi)$ a deformation problem $\mathcal{D}_v$ over $\CO$, in the sense of \cite[Definition 2.2.2]{CHT}, by the condition that,  $\cL_v(\pi')$ is a deformation of $\cL_v(\pi)$ of type $\mathcal{D}_v$ if and only if 
$\Hom_{K_v^+}(\Yu_v,\pi') \neq 0$.  To $\mathcal{D}_v$ is attached a subspace $L_v \subset H^1(W_v,\adj(\overline{\cL}_v(\pi)))$, as in \cite[\S 2.2]{CHT}, and we  assume that 
the problem $\mathcal{D}_v$ is minimal in the sense that
\begin{equation}\label{mini} \dim L_v = \dim H^0(W_v,\adj(\overline{\cL}_v(\pi))).
\end{equation}
Since we have already assumed in Hypothesis \ref{Yuhyp} (2) that $p$ is a banal prime for the local group, some of these hypotheses may be 
redundant.

\section{Perfect complexes and diamond operators}

\subsection{Review of the theorem of Nakajima}\label{naka}
We begin with a simple generalization of a special case of the theorem of Nakajima \cite[Theorem 2]{N} used in \cite{H13}.

\begin{thm}\label{nakajima}    Let $\Gamma$ be a finite abelian $p$-group.  Let  $k$ be a field of characteristic $p$ and let $f:  X \ra Y$ a finite \'etale
Galois covering of projective varieties  over $k$ with Galois group $\Gamma$.  Let $\CF$ be a coherent
sheaf on $Y$.  Let $T$ be a commutative $k$-algebra of endomorphisms of the cohomology complex $R\Gamma(X, f^*(\CF))$ in the derived category of $k[\Gamma]$-modules,
and let $\fm \subset T$ be a maximal ideal with the property that 
$H^i(R\Gamma(X, f^*(\CF))_{\fm}) = 0$ for all indices except $i = i_0$, where the subscript $\fm$ denotes localization at $\fm$.  Then
$H^{i_0}(X,f^*(\CF))_{\fm}$ is a free $k[\Gamma]$-module.  
\end{thm}
\begin{proof}  Theorem 1 of \cite{N} implies that $R\Gamma(X, f^*(\CF))$ can be represented by a finite complex of projective $k[\Gamma]$-modules $C^{\bullet}$.  The localization $C^{\bullet}_{\fm}$ is a direct summand of $C^{\bullet}$ and is therefore also a finite complex of projective $k[\Gamma]$ modules.   It then follows, as in the proof of Theorem 2 of \cite{N}, that 
$H^{i_0}(X,f^*(\CF))_{\fm}$ is a projective $k[\Gamma]$-module, hence free because $k[\Gamma]$ is a local ring.

\end{proof}  

   The following corollary is then proved just as in \cite[Corollary 3.5]{H13}:
\begin{cor}\label{nakajimaCO}  Let $\CO$ be a $p$-adic integer ring with residue field $k$, and let
$f:  X \ra Y$ be a finite \'etale Galois covering of projective $\CO$-schemes with Galois group $\Gamma$, a finite abelian $p$-group.
 Let $T$ be a commutative $\CO$-algebra of endomorphisms of the cohomology complex $R\Gamma(X, f^*(\CF))$ in the derived category of $\CO[\Gamma]$-modules,
and let $\fm \subset T$ be a maximal ideal with the property that 
$$H^i(R\Gamma(X, f^*(\CF))_{\fm}) = 0$$ for all indices except $i = i_0$, where the subscript $\fm$ denotes localization at $\fm$.  
 Assume moreover
that $H^{i_0}(X,f^*(\CF))_{\fm}$ is $\CO$-torsion-free.  Then
$H^{i_0}(X,f^*(\CF))_{\fm}$ is a free $\CO[\Gamma]$-module.
\end{cor}

\subsection{}\label{diamond}.  We apply this to the situation of Proposition \ref{torsionfree}.   Fix a neat level subgroup $K$ and let
 $Q$ be a collection of primes of $F^+$, split in $F/F^+$, at which $K$ is hyperspecial maximal.  Choose
toroidal data $\Sigma_{i,Q}$ as in Theorem \ref{etale}, and write
$$\bbS_{i,Q,\Sigma} = \bbS_{K_{i,Q}(t)}(G,X)_{\Sigma_{i,Q}(t)}, ~~ i = 0, 1, \Delta.$$ 
For $\mu \in X^+(T)$ and $w \in W^x$ we write $\CE_w(\mu)$ for the automorphic vector bundle on $\bbS_{i,Q,\Sigma}$ (any $i$)
attached to the representation of $K_x$ with highest weight $w\cdot \mu$ -- this is the bundle that was denoted
$\CE_w(W)$ in \cite{H13}, where $W$ is the representation of $G$ with highest weight $\mu$.   Let $\CE_w(\mu)^{\rm can}$ denote its canonical extension on $\bbS_{i,Q,\Sigma}$.  Let $i_0(w,\mu) = i_0(\CE_w(\mu))$, in the notation of Proposition \ref{torsionfree}.
We let $T_{S,K}(w,\mu)$ be the algebra denoted $T_{S,K,\CF}$ in \S \ref{Heckealgebras}, for the automorphic
bundle $\CF_{K} = \CE_w(\mu)$.

The following Corollary extends \cite[Proposition 3.7]{H13} to the noncompact case.

\begin{cor}\label{freediamond}  Let $\mu \in X^+(T)$ be a character and let $\fm \subset T_{S,K}(w,\mu)$ be a maximal ideal that
satisfy the hypotheses of Proposition \ref{torsionfree}.   Then for all
$w \in W^{x}$, 
$$H^{i_0(w,\mu)}(\mathbb{S}_{\Delta,Q,\Sigma}, \CE_w(\mu)^{\rm can})_{\fm}$$ is a free 
 $\mathcal{O}\left[\Delta_{Q}\right]$-module
and
$$
H^{i}(\mathbb{S}_{\Delta,Q,\Sigma}, \CE_w(\mu)^{\rm can})_{\fm} =0 \quad \text { if } i \neq i_0(w,\mu).
$$
\end{cor}  

\begin{proof}  In view of Theorem \ref{etale}, we can apply the same  proof as in \cite[Proposition 3.7]{H13}.
\end{proof}

More generally, we introduce a set of type data as in \S \ref{typedata}.   We let $T_{S,K}(w,\mu,\Yu)$ be the analogue of $T_{S,K}(w,\mu)$ for the automorphic
bundle $\CF_{K,\Yu} = \CE_w(\mu)\otimes_\CO \Lambda_{\Yu}$.  Then we have

\begin{cor}\label{freediamondYu}  Let $\mu \in X^+(T)$ be a character and let $\fm \subset T_{S,K}(w,\mu)$ be a maximal ideal that
satisfy the hypotheses of Proposition \ref{torsionfree}.   Then for all
$w \in W^{x}$, 
$$H^{i_0(w,\mu)}(S_{\Delta}(Q), \CE_w(\mu)^{\rm can}\otimes_\CO \Lambda_{\Yu})_{\fm}$$ is a free 
 $\mathcal{O}\left[\Delta_{Q}\right]$-module
and
$$
H^{i}(\mathbb{S}_{\Delta}(Q), \CE_w(\mu)^{\rm can}\otimes_\CO \Lambda_{\Yu})_{\fm,\Yu} =0 \quad \text { if } i \neq i_0(w,\mu).
$$
\end{cor}  
The proof is the same as in the case without type data; we only need to add that the functor $\Hom_{K_S^+}(\Lambda_{\Yu},\bullet)$
is exact because we have assumed that $p$ is a banal prime for all places in $S^+$.

In the applications we impose the conditions of \S \ref{minimality} on our choice of type data to ensure that the local deformation conditions are minimal.

\section{Application of the Taylor-Wiles method}  We follow the axiomatic treatment of the Taylor-Wiles method
described in \S 4 of \cite{H13}.     Fix $\mu$ and $w \in W^x$ as in Corollary \ref{freediamond}, and let $i_0 = i_0(w,\mu)$
for the duration of this section.  Fix a neat level subgroup $K$, define $S = S(K)$ as above, and consider sets $Q$ and a maximal ideal
$\fm$ as in \S \ref{diamond}.

Let
$$H_\emptyset = H^{i_0}(\bbS_{K,\Sigma},\CE_w(\mu)^{\rm can}),
H_{0,Q} = H^{i_0}(\bbS_{0,Q,\Sigma},\CE_w(\mu)^{\rm can}),$$ 
$$H_{\Delta,Q} = H^{i_0}(\bbS_{\Delta,Q,\Sigma},\CE_w(\mu)^{\rm can}).$$

Note that by Proposition~ \ref{torsionfree}, these are torsion-free $\CO$-modules.  The Hecke algebras
$\TT_{\emptyset,w}, \TT_{0,Q,w}, \TT_{\Delta,Q,w}$ are defined as in \cite[p. 137]{H13} with respect to these 
modules as the images of the relevant Hecke algebras $\TT_{S,K}(w,\mu)$, initially defined as endomorphisms
of the perfect complex, on the cohomology modules.  Given a homomorphism 
$$\nu:  \TT_{\emptyset,w} \ra \Qpb$$ 
we let $\rho_\nu$ be the corresponding Galois representation, and let $R_{\bar{\rho}_\nu,\emptyset}$ denote
the corresponding deformation ring when (as we assume below) $\bar{\rho}_{\nu}$ is absolutely irreducible.

We have the following analogue of Theorem 6.8 of \cite{H13}:
\begin{thm}\label{TWtheorem}
Assume $\mu$ satisfies the inequalities of Proposition~ \ref{torsionfree}.  
Let $\nu:  \TT_{\emptyset,w}(\mu) \ra \Qpb$ be a non-trivial homomorphism
and assume the corresponding Galois representation $\rho_{\nu}$ satisfies conditions
\begin{enumerate}
\item For all $v \in S(F^+)$, $\bar{\rho_{\nu}}$ is unramified at $v$ and 
\[
H^0( \Gal(\bar{F_v}/F_v),  (\textrm{\em ad} ~\bar{\rho_{\nu}})(1))=(0);
\]
%\item For all $q \in S(G)$, $\bar{\rho_{\nu}}$ is unramified for all $v \in D_-(q)$ and 
%\[
%H^0( \Gal(\bar{F_v}/F_v),  (\textrm{\em ad} ~\bar{\rho_{\nu}})(1))=(0);
%\] for such $v$.
\item The Fontaine-Laffaille condition at primes above $p$.
\end{enumerate}

We also assume the residual representation $\bar{\rho}_{\nu}$ to be absolutely irreducible\footnote{needed for the construction of the sets $Q_M$ in Prop. 4.4 by Thorne.} and {\it adequate}
in the sense of \cite{Th}.  Let $\fm = \fm(\bar{\rho}_\nu) \subset \TT_{\emptyset,w}(\mu)$ be the corresponding
maximal ideal, and let $H_{\bar{\rho}_\nu,\emptyset}$ and $\TT_{\bar{\rho}_\nu,\emptyset}$ denote the localizations at $\fm$.
Suppose that the deformation conditions are \it{minimally ramified} at all primes $v\in S$ not dividing $p$. In particular, for $v \in S$, we assume the deformation condition in $\CS$ at $v$
is unrestricted and minimal (cf. \cite[\S 6.9]{H13}). 

Then the classifying map
$$\phi_{\bar{\rho}_\nu,\emptyset}:  R_{\bar{\rho}_\nu,\emptyset} \ra \TT_{\bar{\rho}_\nu,\emptyset}$$
is an isomorphism, and 
$H_{\bar{\rho}_\nu,\emptyset}$ is a free module over $\TT_{\bar{\rho}_\nu,\emptyset}$, which is a local complete intersection.

More generally, suppose we replace $\CF_K$ by $\CF_{K,\Yu}$, for type data as in \S \ref{typedata}. Let  $\TT_{\emptyset,w,\Yu}(\mu)$ be the corresponding Hecke algebra, let $\nu:  \TT_{\emptyset,w,\Yu}(\mu) \ra \Qpb$
be a homomorphism satisfying the conditions above, plus the minimality condition of \eqref{mini} at places in $S$.
Define $\TT_{\bar{\rho}_\nu,\emptyset}$ and $R_{\bar{\rho}_\nu,\emptyset}$ with respect to this homomorphism $\nu$.  Then the conclusions of the theorem stated above
remain true.
\end{thm}

\begin{rmk} The reader can check that, since $p$ is odd, the condition on $v \in S(F^+)$ is not strictly necessary; this is the point of 
\cite[Lemma 2.3.3]{CHT}, whose application is buried at the end of the proof of Lemma 2.3.4 of the same paper.  In order to avoid distraction 
we add this as an extra condition.  %The condition on $q \in S(G)$ is introduced similarly in order to avoid 
\end{rmk}

%{\color{red}  Have to show that the localization is reduced here.}

%\begin{rmk}\label{neat1}\end{rmk}

\begin{proof}
We show that the relevant hypotheses in  \S 4 of \cite{H13} are satisfied with the necessary adjustment to account for the assumption
that our residual representations are adequate. The treatment of \cite{H13} follows the Diamond-Fujiwara extension to the Taylor-Wiles method as applied in~\cite{CHT}, while we use a modification as in Theorem 3.6.1 in \cite{BLGG}.

The Galois Hypothesis 4.3 follows as in the proofs of Corollaries 5.3 and 5.4 in \cite{H13} in light of our Theorem~\ref{Galoisreps} (ii). %Note that $\TT_{\emptyset,w}(\mu)$ may not be reduced; for that reason in all proofs we use the \emph{reduced localizations}. 
Hypothesis 4.4.8 is verified in Lemma~6.1 therein. It remains to check the correct analog of Hypothesis 4.4.9 for adequate residual representations.

The conditions (4.4.9.1) and (4.4.9.2) in Part (a) of Hypothesis 4.4.9 follow from the constructions of sets of Taylor-Wiles $Q_M$ satisfying the conditions (6.6.1) and (6.6.2) from the beginning of \S6.6 in \cite{H13}; such a set $Q_M$ is shown to exist in Proposition~4.4, \cite{Th}. The hypotheses of \S~4.4.5 and \S~4.4.7 as well as Hypothesis 4.4.9 (b) are verified for the set $Q_M$ in the proof of Theorem~6.8 in \cite{Th}. In the same proof one can find the verification of the local condition in \S~4.4.4. as a consequence of Proposition~5.12 in \cite{Th}. The modified version of Condition 4.4.9.3, now concerning nonvanishing of trace of Frobenius along projection, is embedded in the construction of $Q_M$ in Proposition~4.4. The claim now follows since under our assumptions (i), (ii) and (iii) in the proof of Theorem~6.8~\cite{Th} hold.

Finally, the modifications necessary for incorporating the type data are the same as in the proof of \cite[Theorem 3.5.1]{CHT}, which specifically comes down to the use
of the minimality hypothesis in \cite[Proposition 2.5.9]{CHT}.
\end{proof}

%\begin{rmk}[Erratum for \cite{H13}]  The condition (6.3) in 

%\end{rmk}
\section{Consequences for cohomology modules over the Hecke algebra}\label{modules}

\subsection{Topological and de Rham cohomology}  
Theorem \ref{TWtheorem} asserts that the coherent cohomology module $H^{i_0}(\bbS_{K,\Sigma},\CE_w(\mu)^{\rm can})$ localized at the maximal ideal $\fm = \fm(\bar{\rho}_\nu)$ is free over the localized Hecke algebra $\TT_{\bar{\rho}_\nu,\emptyset}$, provided  $p$ and the weight $\mu$ satisfy the stated inequalities.
As in \cite[\S 7]{H13}, this implies that analogous localized modules in $p$-adic 
\'etale and de Rham cohomology are also free over $\TT_{\bar{\rho}_\nu,\emptyset}$.  The proofs are identical, so we merely state the results.

Let now $W = W_\mu$ be the finite-dimensional irreducible representation of $G$ with highest weight $\mu$.  As in \cite[(7.1.2)]{H13} we always assume
$\mu$ to be $p$-small.  Let the $p$-adic place $v$, the $p$-adic integer ring $\CO$, the hyperspecial maximal compact subgroup $K_p \subset G(\Qp)$,
and the $K_p$-invariant lattice $W_\mu(\CO) \subset W(Frac(\CO))$ be as in \cite[\S 7.1]{H13}.  Let $\tilde{W}_{\mu, B}(\mathcal{O})$
denote the topological local system in $\CO$-modules on $S_{K}(G, X)(\mathbb{C})$, as in \cite[\S 7.3]{H13}.  Topological cohomology is computed by
the complex 
$$R\Gamma(S_{K}(G, X)(\mathbb{C}), \tilde{W}_{\mu, B}(\mathcal{O})).$$  We let 
$$\TT_{\emptyset,B}(\mu) \subset End \left(R\Gamma(S_{K}(G, X)(\mathbb{C}), \tilde{W}_{\mu, B}(\mathcal{O})\right)$$
denote the $\CO$-subalgebra
generated by the topological Hecke correspondences $T(t)$ 
introduced in \S \ref{Heckealgebras}.  

The de Rham version $\TT_{\emptyset,dR}(\mu)$ of $\TT_{\emptyset,B}(\mu)$ requires a bit more work to define.  Let $\CE_{G,\CO}$ be the exact monoidal functor 
defined in \cite[Lemma 1.20]{LS1} from the tensor category of representations of the group scheme $G_\Zp$ over $\CO$ to the category of locally free coherent sheaves
on $\bbS_K$ with integrable connection and compatible $G(\ad_f^p)$ action covering the natural action on the family $\bbS_{K\times K^p}$.  Let $W_{\mu,\CO}$ denote
the representation of $G_\Zp$ introduced in \cite[\S 7.1]{H13}, so that the $\CO$-lattice $W_\mu(\CO)$ is just $W_{\mu,\CO}(\CO)$.  We let
\begin{equation}\label{RdR}
R\Gamma_{dR,log}(\bbS_K,W_\CO) = R\Gamma(\bbS_{K,\Sigma},\Omega_{\bbS_{K,\Sigma},log}^\bullet \otimes \CE_{G,\CO}(W_{\mu,\CO})^{can});
\end{equation}
see \cite[\S 4.3]{LS2}, or \cite[\S 7.4]{LS2} for a reformulation in terms of dual BGG complexes.   We define 
 $\TT_{\emptyset,dR}(\mu) \subset End\left(R\Gamma_{dR,log}(\bbS_{K}, W_\CO))\right)$
 again to be the  $\CO$-subalgebra generated by the algebraic Hecke correspondences $T(t)$.

We let $\nu_B:  \TT_{\emptyset,B}(\mu) \ra \Qpb$ (resp. $\nu_{dR}:  \TT_{\emptyset,dR}(\mu) \ra \Qpb$) be a non-trivial homomorphism
satisfying the hypotheses of Theorem~ \ref{TWtheorem}.  Define $\fm_B = \fm(\bar{\rho}_\nu) \subset \TT_{\emptyset,B}(\mu)$ (resp.  
$\fm_{dR} = \fm(\bar{\rho}_\nu) \subset \TT_{\emptyset,dR}(\mu)$ as in the statement of that theorem, and let
$\TT_{\bar{\rho}_\nu,\emptyset,B}(\mu)$, $\TT_{\bar{\rho}_\nu,\emptyset,dR}(\mu)$ denote the respective localizations.   We let
$$R\Gamma_{dR,log}(\bbS_{K}, W_\CO))_{\bar{\rho}_\nu}, ~~ R\Gamma(S_{K}(G, X)(\mathbb{C}), \tilde{W}_{\mu, B}(\mathcal{O}))_{\bar{\rho}_\nu}$$
denote the  localizations of the cohomology complexes at the respective ideals.

\begin{hyp}\label{runninghypotheses}  In what follows, we assume
\begin{itemize}
\item[(a)] $\mu$ satisfies the inequalities of Proposition~ \ref{torsionfree};
\item[(b)] $\bar{\rho}_\nu$ satisfies the hypotheses of Theorem \ref{TWtheorem}.
\end{itemize}
\end{hyp}

Recall that $d_V = \dim S_K(G,X)$.  

\begin{cor}\label{free}  Under Hypotheses \ref{runninghypotheses}  we have
(i)
$$R\Gamma_{dR,log}(\bbS_{K}, W_\CO)_{\bar{\rho}_\nu}$$ 
is concentrated in degree $d_V$, and
$$H^{d_V}_{dR,log}(\bbS_{K}, W_\CO)_{\bar{\rho}_\nu} := H^{d_V}(R\Gamma_{dR,log}(\bbS_{K}, W_\CO)_{\bar{\rho}_\nu})$$
is a free $\CO$-module of finite rank.  Moreover, there is a natural decreasing (Hodge) filtration
$F^\bullet H^{d_V}_{dR,log}(\bbS_{K}, W_\CO))_{\bar{\rho}_\nu}$ by $\CO$-direct summands satisfying the analogue
of \cite[Theorem 7.2.2]{H13}.

(ii) Assume in addition the inequality $|\mu_{comp}| \leq p - 2$, as in the statement of \cite[Theorem 7.3.3]{H13}.  Then 
$$R\Gamma(S_{K}(G, X)(\mathbb{C}), \tilde{W}_{\mu, B}(\mathcal{O}))_{\bar{\rho}_\nu}$$
is concentrated in degree $d_V$.  Moreover  $H^{d_V}(S_{K}(G, X)(\mathbb{C}), \tilde{W}_{\mu, B}(\mathcal{O}))_{\bar{\rho}_\nu}$
is a free $\CO$-module of finite rank (the same as in (ii)).
\end{cor}
\begin{proof}  This follows from the earlier results in the same way as in \cite{H13}, except that here we need to localize at the non-Eisenstein maximal
ideal $\fm(\bar{\rho}_\nu)$.
\end{proof}

The next two theorems are then proved exactly as in \cite[\S 7.2, \S 7.3]{H13}.
\begin{thm}\label{toptheorem}  Assume Hypotheses \ref{runninghypotheses}   and  the inequality $|\mu_{comp}| \leq p - 2$, as in the statement of \cite[Theorem 7.3.3]{H13}.  
Then there is an isomorphism
$$\phi_{\bar{\rho}_\nu,\emptyset,B}:  R_{\bar{\rho}_\nu,\emptyset} \ra  \TT_{\bar{\rho}_\nu,\emptyset,B}(\mu)$$
of local complete intersections, and $H^{d_V}(S_{K}(G, X)(\mathbb{C}), \tilde{W}_{\mu, B}(\mathcal{O}))_{\bar{\rho}_\nu}$ is a free
module of finite rank over $\TT_{\bar{\rho}_\nu,\emptyset,B}$.
\end{thm}

We remind the reader that the inequality $|\mu_{comp}| \leq p - 2$ is needed in order to apply a comparison theorem for integral $p$-adic cohomology.

\begin{thm}\label{dRtheorem} Assume $\mu$ satisfies the inequalities of Proposition~ \ref{torsionfree}.   
Then there is an isomorphism
$$\phi_{\bar{\rho}_\nu,\emptyset,dR}:  R_{\bar{\rho}_\nu,\emptyset} \ra  \TT_{\bar{\rho}_\nu,\emptyset,dR}(\mu)$$
of local complete intersections, and $H^{d_V}_{dR,log}(\bbS_{K}, W_\CO))_{\bar{\rho}_\nu}$ is a free
module of finite rank over $\TT_{\bar{\rho}_\nu,\emptyset,dR}(\mu)$.
\end{thm}

\subsection{Assumptions}\label{as}
Let $\TT$ be a finite free reduced local $\CO$-algebra, $M$ a finite free $\CO$-module, with $\TT$ action
defined by a faithful homomorphism $\imath:  \TT \hookrightarrow End_\CO(M)$.  
Then $\TT_\QQ = \TT\otimes_\Zp \Qp$ is a product of finite extensions of $\Qp$ that acts faithfully on
$M_\QQ = M\otimes_\Zp \Qp$.  We assume $M_\QQ$ is free over $\TT_\QQ$ of rank $m$.  For $0 \leq i \leq m$ we consider
$$\wedge^i M_\QQ = \wedge^i_{\TT_\QQ} M_\QQ;$$
this is naturally a free $\TT_\QQ$-module of rank $\frac{m!}{i!\cdot (m-i)!}$.

In what follows we assume $p > m$, in order to avoid ambiguity regarding the presence or not of factorials in the definition of exterior powers over $\CO$-algebras.  

\subsubsection{Multiplicity one}  With a bit more work we can show that we can take $m = n$ when $M =  H^{n-1}_{dR,log}(\bbS_{K}, W_\CO)_{\bar{\rho}_\nu}$,
as in the next section.  It would be convenient to assume that cuspidal automorphic
representations of $G$ occur that are spherical outside primes that split in $F/F^+$  satisfy strong multiplicity one:  they are determined uniquely, as subspaces
of the space of cusp forms, by their local components at all unramified places.  Because the ramification is limited to places where $G$ is isomorphic (up to the similitude factor)
to $GL(n)$, this follows -- almost -- from the multiplicity one theorem for $L$-packets of unitary groups proved in \cite{KMSW} (conditionally on unpublished results of Arthur).

The ``almost" refers to the structure of unramified $L$-packets at primes that ramify in $F/F^+$.  There the $L$-packets can have several members, distinguished by the choice of local maximal compact $K_q$; so we recover multiplicity one, though not necessarily strong multiplicity one.  This, together with the properties of Bushnell-Kutzko types recalled in Definition \ref{semisimpletype}, is enough to obtain $m = n$ in the situation considered after
Proposition \ref{ext}.

\subsection{Exterior powers of modules over Hecke algebras}
{~}

\begin{prop}\label{ext}  With notation as in \S \ref{as}, assume $M$ is free over $\TT$.  Let $M_i \subset \wedge^i M_\QQ$ be a free $\TT$-submodule
of rank $\frac{m!}{i!\cdot (m-i)!}$.  Then $M_i$ and $\wedge^i_{\TT} M$
are isomorphic as $\TT$-modules.
\end{prop}

\begin{proof}  
Since $M$ is a free $\TT$-module of rank $m$, then $\wedge^i_{\TT} M$ is also a free $\TT$-module of rank ${m \choose i}$ (see Theorem~4.2 in \cite{KC}). Now, $M_i$ and $\wedge^i_{\TT} M$ are finite free modules of the same rank, hence isomorphic.
\end{proof}

We apply this when $M = H^{n-1}_{dR,log}(\bbS_{K}, W_\CO)_{\bar{\rho}_\nu}$  in the notation of Corollary \ref{free},
with $\mathbb{S}_K$ an integral model of the Shimura variety $Sh(V_\sigma)$ attached to an $n$-dimensional hermitian space $V_\sigma$ over $F$ with signature $(n-1,1)$ at one place  $\sigma \in \tilde{S}_\infty$, 
and definite at the remaining places, and with $\TT = \TT_{\bar{\rho}_\nu,\emptyset,dR}(\mu)$ the corresponding localized Hecke algebra.  (If $F$ is not imaginary quadratic the
subscript $log$ is irrelevant.)  We also write $M = M_\sigma$ and $\TT = \TT_\sigma$ and allow $\sigma$ to vary over $\tilde{S}_\infty$; we write $G_\sigma$ for the
similitude group of $V_\sigma$ and $(G_\sigma,X_\sigma)$ for the corresponding Shimura datum.  

The algebra $\TT_\sigma$ is the localization of the full Hecke algebra $ \TT_{\emptyset,dR}(\mu)$ at the maximal ideal containing the kernel of  the homomorphism $\nu_{dR}$;
we use the same notation to denote the restriction $\nu_{dR}:  \TT_\sigma \ra \Qpb$ to the localized Hecke algebra.  This homomorphism corresponds to the action of
$\TT_\sigma$ on the $K$-fixed vectors in an automorphic representation $\pi = \pi_\nu$ of $G_\sigma(\ad)$ that is realized in  $H^{n-1}_{dR,log}(S_K(G_\sigma,X_\sigma), W_\mu)$.

Now let $V$ be any $n$-dimensional (non-degenerate) hermitian space over $F$, $(G_V,X_V)$ the corresponding Shimura datum, $K_V \subset G_V(\af)$ a neat level subgroup that is hyperspecial maximal at $p$,
$Sh_{K_V}(V)$ the corresponding Shimura variety, $\bbS_{K_V,V}$ a smooth model over $\CO$ (we may have to take $\CO$ sufficiently large to include the integer rings in all
relevant reflex fields).   We let $\TT_{\emptyset,dR,V}(\mu)$ denote the $\CO$-algebra generated by  Hecke operators  at places split in $F/F^+$ acting on
$H^{d_V}_{dR,log}(\bbS_{K_V,V}, W_\CO)$, with $W = W_\mu$; the index $K_V$ is omitted from the notation for the Hecke algebra.  We can discard the Hecke operators
at a finite number of places without changing $\TT_{\emptyset,dR,V}(\mu)$, and thus we can define $\nu:  \TT_{\emptyset,dR,V}(\mu) \ra \Qpb$ and the localization
$$\TT_V := \TT_{\bar{\rho}_\nu,\emptyset,dR,V}(\mu)$$
 for any $V$ such that the representation $\pi_\nu$ transfers to $G_V(\ad)$.  

The congruence ideal $C(\nu,\TT_V) = C(\pi_\nu,\TT_V) \subset \TT_V$ is defined as in \cite[Definition 6.7.2]{EHLS}.  More precisely, it is the annihilator of
\[
\TT_V/(\TT_V[\pi_\nu]+\TT_V[\pi_\nu]^{\bot}).
\]
Here $\bar{\pi}_\nu$ is the dual of $\pi_\nu$,  $\TT_V[\pi_\nu]$ is the localization of $\TT_V$ at the prime ideal that is the kernel of the action of
the Hecke algebra on the $\pi_\nu$-isotypic subspace in $\TT_V\otimes \QQ$, and $\TT_V[\pi_\nu]^{\bot}$ is the intersection of $\TT_V$
with the orthogonal complement of $\TT_V[\bar{\pi}_\nu]$ in $\TT_V\otimes \QQ$  with respect to Poincar\'e duality.  

\begin{prop}\label{congruenceV}  Let $V$ and $V'$ be two hermitian spaces of dimension $n$ such that the representation $\pi_\nu$ transfers to $G_V(\ad)$
and $G_{V'}(\ad)$.   
 Then we can identify the Hecke algebras $\TT_V \isoarrow \TT_{V'}$
 in such a way that the Hecke operators at places $v$ such that  both $G_{V}$ and $G_{V'}$ are split at $v$ correspond.  Moreover, suppose   $\mu$ and $p$ satisfy the inequalities of Proposition~ \ref{torsionfree} for both $d_V$ and $d_{V'}$.    Then with respect to the identification of $\TT_V$ with $\TT_{V'}$, 
the congruence ideals $C(\nu,\TT_V)$ and $C(\nu,\TT_{V'})$ coincide.
\end{prop}

\begin{proof} The first assertion is clear.  The claim about the congruence ideals is then a consequence of Theorem \ref{dRtheorem}, which identifies both $\TT_V$ and $\TT_{V'}$ with the same $\CO$-algebra $R_{\bar{\rho}_\nu,\emptyset}$.
\end{proof}

In the same way, we see that the congruence ideals defined with respect to the integral structure on de Rham cohomology coincide with the ideals defined with respect to the
graded pieces with respect to the Hodge filtration, which are defined with respect to the integral structure on coherent cohomology.  We leave it to the reader to formulate
the statement.

\begin{rmk}  It suffices to say that all the results of this section remain valid when the coefficients $W_{\mu,\CO}$ are replaced by $W_{\mu,\CO}\otimes \Lambda_{\Yu}$ for type data satisfying the minimality conditions, and
the (de Rham or $p$-adic \'etale) cohomology modules $\bullet = R\Gamma_{dR,log}(\bbS_K,W_\CO), \dots$ are replaced by $\Hom_{K_S^+}(\Lambda_{\Yu},\bullet)$.
\end{rmk}

\section{Ordinary modular forms}\label{ordinarymf}

The article \cite{EHLS} constructs $p$-adic $L$-functions as elements of Hida's ordinary Hecke algebra when the latter, localized at
an appropriate maximal ideal, is known to satisfy a {\it Gorenstein hypothesis}  \cite[\S 6.7.8, \S 7.3.2]{EHLS}.    Here we show how to derive
this hypothesis from Theorem \ref{TWtheorem}, when $p$ is sufficiently large and the maximal ideal is as in the statement of that theorem.  The
first results of this type are due to Hida \cite{Hi86} and Tilouine \cite{T87}; our method is essentially the same as Tilouine's.

We make use of the notation of \cite{EHLS} without comment.  Thus $\pi$ is an anti-holomorphic cuspidal automorphic representation of the unitary similitude
group $G$ (denoted $G_1$ in {\it loc. cit.}) and $\TT = \TT_\pi$ is the localization of Hida's ordinary Hecke algebra at the corresponding maximal ideal $\fm_\pi$.  This
ideal determines a connected component of weight space, which includes  the weight $\kappa(\pi)$ of the holomorphic modular form corresponding to $\pi$ (more conventionally, to the contragredient of $\pi$, a holomorphic cuspidal automorphic representation), which we are allowed to vary.  Theorem \ref{TWtheorem} can be applied provided $\pi$ can be chosen of a weight that satisfies the inequalities in Proposition \ref{torsionfree} as well as those of Hida's Control Theorem, cited as \cite[Theorem 7.2.1]{EHLS}.  We make these conditions more precise in the following Theorem.   We say two weights $\mu, \mu'$ are congruent modulo $p$ if they define the same characters
on the (finite group of) points of the special fiber of the torus $T$ over the residue field $k$ of $\CO$.

\begin{thm}\label{Gorenstein}   Let $\kappa(\pi)$ be a weight such that  the $\CO$-module $M_{\kappa(\pi)}(K)$ of holomorphic modular forms of weight $\kappa(\pi)$, in
the notation of \cite{EHLS}, is the module $H^{i_0(w,\mu)}(\bbS_{K,\Sigma},\CE_w(\mu)^{\rm can})$, in the notation of Corollary \ref{freediamond}, with $i_0(w,\mu) = 0$.  
Suppose $\mu$ is congruent modulo $p$ to a character $\mu' $  of $T$  that satisfies the regularity condition (ii) of Proposition \ref{torsionfree}.   Suppose the residual Galois representation $\bar{\rho}_\pi$ attached to the maximal ideal 
$\fm_\pi$ satisfies the hypotheses of Theorem \ref{TWtheorem}.
Then the localized ordinary Hecke algebra $\TT = \TT_\pi$ is a local complete intersection and satisfies the Gorenstein Hypothesis 7.3.2 of \cite{EHLS}.
\end{thm}

\begin{proof}  The notation $\bar{\rho}_\pi$ corresponds to the notation $\bar{\rho}_\nu$ in the statement of Theorem \ref{TWtheorem}.
Applying that theorem, the hypotheses imply that the localized module $S_{\kappa(\pi)}(K)_\pi$ is free over the localized Hecke algebra, which is a local complete intersection.
We claim that it suffices to show that Hida's Control Theorem, as stated in \cite[Theorem 7.1]{Hi02}, applies to $S_{\kappa(\pi)}(K)_\pi$.  Indeed, the Hecke algebra
$\TT_{\bar{\rho}_\nu,\emptyset}$ in weight $\kappa(\pi)$ is the quotient of $\TT_\pi$ by the regular sequence corresponding to the weight $\kappa(\pi)$.  Moreover, 
we are concerned with coherent cohomology in degree $0$, so by Koecher's principle the reduced localization hypothesis in Theorem \ref{TWtheorem} 
is superfluous.

Now it remains to verify that Hida's Control Theorem does apply to weight $\kappa(\pi)$ when $\mu'$ satisfies the regularity condition (ii) of Proposition \ref{torsionfree}.
In fact, Boxer and Pilloni have proved a classicality theorem for overconvergent modular forms of small slope \cite[Theorem 1.0.15 (3)]{BP} that asserts that
overconvergent  modular forms of weight $\mu$ on Shimura varieties of abelian type are classical if they satisfy a small slope condition and if $\mu'$  satisfies the 
regularity condition (ii) of Proposition \ref{torsionfree}.  Since ordinary modular forms automatically satisfy the small slope condition (see \cite[Remark 1.0.6]{BP}),
this completes the proof.

\end{proof}

\begin{rmk}\label{generalmu}  It is important to note that Theorem \ref{Gorenstein} is really a statement about $\bar{\rho}_\pi$ and the prime $p$, and not about $\mu$.  If $p$ is not too small relative to the group $G$, then every $\mu$ is congruent to a $\mu'$ satisfying the regularity condition.  
\end{rmk}

\begin{rmk}  Once again, it suffices to say that all the results of this section remain valid when the coefficients $\CE_w(\mu)^{\rm can}$ are replaced by $\CE_w(\mu)^{\rm can}\otimes \Lambda_{\Yu}$ for type data satisfying the minimality conditions and the cohomology modules are modified accordingly.  We leave the details to the reader.
\end{rmk}

\section{Proof of Theorem \ref{Galoisreps} (iii) and (iv) }\label{proofofreps}

The proof that the Galois representations attached to the boundary cohomology group 
$R\Gamma^{\partial}(\CF_K) := R\Gamma^{\partial}(\bbS_K(G,X),\CF_K)$,
or to its generalization $$R\Gamma^{\partial}(\bbS_K(G,X),\CF_{K,\Yu}),$$
are reducible is based on a lengthy analysis of the toroidal boundary, to show that the
Galois representations attached to any piece of the complex that computes the cohomology of the toroidal 
boundary breaks up as a sum of representations attached to the cohomology of smaller groups.  This is a geometric
argument and it is identical with or without the introduction of a $K$-type indicated by the subscript $\Yu$.  Thus we
will write down the proof, which is long enough as it is, in the case where the type data are trivial.

 In what follows, $R$ denotes a standard rational proper parabolic
subgroup of $G$.  We let $r(R)$ denote the parabolic rank of $R$; thus $r(R) = 1$ if and only if
$R$ is a maximal proper parabolic.   Let $L_R$ denote the Levi quotient of $R$.  It is a connected reductive
group that admits a factorization
$$L_R = G_{h,R} \cdot G_{\ell,R}$$
where $G_{h,R} = G(V_R)$ is the similitude group of a hermitian vector space $V_R$ over $F$, of dimension $n - 2m(R)$, for
some $1 \leq m(R) \leq \frac{n}{2}$, and 
$$G_{\ell,R} = \prod_{i = 1}^{r(R)} GL(m_i(R))_F, \sum_i m_i(R) = m(R).$$
The factorization is defined by choosing a flag 
$$0 = A_0 \subset A_1 \subset A_2 \subset \dots \subset A_{r(R)}$$
of totally isotropic subspaces of $V$, all assumed to be defined over $\CO$.  Then $GL(m_i(R))$ is identified with the group 
scheme $GL(A_i/A_{i-1})$. over $Spec(\CO)$.  In particular, there is a surjective homomorphism
\begin{equation}\label{lin}  \ell_R:  L_R \ra G_{\ell,R} 
\end{equation}
whose kernel is isomorphic to the unitary similitude group $G(V_R)$, where $V_R$ is the quotient of $V/A_{r(R)}$ by the null space of the induced
hermitian form.  

There is a Shimura datum $(G_{h,R},X(R))$ attached to $R$, of the same kind
as the pair $(G,X)$ with which we began.
As noted in the Appendix, this is {\it not quite} the boundary Shimura datum that corresponds to $R$ in the theory of \cite{P} or
\cite{HZ94,HZ01}.  In those references the group $G_{h,R}$ contains an additional split torus, here relegated to $G_{\ell,R}$, that
is needed in order to define the correct mixed Hodge structure on the boundary cohomology.  For the purposes of this
paper this precision will not be necessary.   The absence of this extra factor is compensated by the Tate twists in 
\eqref{Tatetwist}.\footnote{Pink's theory of mixed Shimura data $(G,X)$ allows for the possibility that $X$ is a union of a finite set of copies of 
hermitian symmetric domains in a flag variety attached to $G$.  Taking this into account requires adjustments that
are discussed in \cite[\S 1.1.7]{HZ01}, and that will be irrelevant for our purposes.} We note for future reference that
\begin{equation}\label{PofR}
(G_{h,R},X(R)) = (G_{h,P(R)},X(P(R))), 
\end{equation}
where, as usual, $P(R)$ is the maximal parabolic subgroup to which $R$ is subordinate, as in \S \ref{parab}.

Let $K_{h,R} \subset G_{h,R}(\RR)$ be the stabilizer of a point $h_R \in X(R)$ -- we may as well assume $h_R$
to be a CM point, though this makes no difference --  and let $K_{\ell,R}$ be a maximal connected subgroup of $G_{\ell,R}(\RR)$ that is compact modulo the center.   Define $K_{R,\infty} = K_{h,R}\times K_{\ell,R}$.  
For any compact open subgroup $K_R \subset L_R(\af)$, we let ${}_{K_R}Y(R)$ denote the locally symmetric space
$L_R(\QQ)\backslash L_R(\ad)/K_{R,\infty}\times K_R.$  We always assume $K_R$ to be {\it neat}, so that 
${}_{K_R}Y(R)$ is smooth and its Borel-Serre compactification is a smooth manifold with corners.

We always assume $K_R$ to have the property that $K_{h,R}^f := K_R\cap G_{h,R}(\af)$ contains a hyperspecial maximal
compact subgroup of $G_{h,R}(\Qp)$.  Let
${}_{K_R}Y_\ell(R)$ denote the locally symmetric space
$$G_{\ell,R}(\QQ)\backslash G_{\ell,R}(\ad)/K_{\ell,R}\times \ell_R(K^f_R),$$
with $\ell_R$ as in \eqref{lin}.  We let
$\bbS_{K_R}(G_{h,R},X(R))$ denote the smooth integral model over $\CO$ of the Shimura variety for $(G_{h,R},X(R))$
of level $K_{h,R}^f$.  

The locally symmetric space ${}_{K_R}Y_\ell(R)$ has a family of local coefficient systems attached to finite-dimensional (algebraic) representations of 
$G_{\ell,R}$, with coefficients in any $\CO$-algebra.   Letting $r_\CW:  G_{\ell,R} \ra \Aut(\CW)$ be such a representation, taken with coefficients in the
$\CO$-algebra $A$, 
we denote by $\tilde{\CW}$ the corresponding
local system with coefficients in $A$:
$$\tilde{\CW} = G_{\ell,R}(\QQ)\backslash G_{\ell,R}(\ad) \times \CW(A)/K_{\ell,R}\times K_\ell(R).$$

Fix $K_R$ unramified at $p$ and let $t\in G_{\ell,R}(\af)$. Let $S=S(K_R)$ temporarily denote the set of ramified places for $K_R$. If $t =(t_w)$ with $t_w\in G_w$, set $K_{R,w}(t) = K_{R,w}\cap tK_{R,w}t^{-1}$, and consider 
\[
K_R(t) \coloneqq \prod_w K_{R,w}(t).
\]
Since $\tilde{\CW}$ is algebraic, it is defined over a finite extension $E(\tilde{\CW})$, and under our assumption admits a model over any $\CO$-algebra $A$. For an open subgroup $K'_R\subseteq K_R$, we have a natural finite covering map
\begin{equation*}
    \begin{split}
        \pi^R_{K_R,K'_R}: {}_{K'_R}Y_\ell(R)\to {}_{K_R}Y_\ell(R).
    \end{split}
\end{equation*}
Let $R\Gamma({}_{K_R}Y_\ell(R), \tilde{\CW})$ be a complex in $D^b(\CO)$ computing $H^*({}_{K_R}Y_\ell(R), \tilde{\CW}).$ The map $\pi^R_{K_R,K'_R}$ induces canonical pullback and direct image functors in $D^b(\CO),$ while the isomorphism $*t^{-1}$ induces quasiisomorphism in $D^b(\CO)$ via pullback. Denote by $T_{S,K(R),\CW}(t) \in \Endo(R\Gamma({}_{K_R}Y_\ell(R),\tilde{\CW}_{K_R}))$ the element given by the composition
\begin{multline}\label{hecke-linear}
R\Gamma({}_{K_R}Y_\ell(R),\tilde{\CW}_{K_R}) \overset{\pi_{K_R,K_R(t)}^*}{\longrightarrow} R\Gamma({}_{K_R(t)}Y_\ell(R),\tilde{\CW}_{K_R(t)}) \\ \overset{\pi_{t^{-1}K_Rt,K_R(t),*}}{\longrightarrow} R\Gamma({}_{t^{-1}K_Rt}Y_\ell(R),\tilde{\CW}_{t^{-1}K_Rt})
\overset{*t^{-1}}{\longrightarrow} R\Gamma({}_{K_R}Y_\ell(R),\tilde{\CW}_{K_R}).
\end{multline}

Denote by $T_{S,K(R),\CW}$ the algebra $\Endo(R\Gamma({}_{K_R}Y_\ell(R),\tilde{\CW}_{K_R}))$ spanned by the above operators with $t\in G_q$ with $w\not\in S(K)\cup \{p\}$ and with the further stipulation that all $w|q$ in $F^+$ split in $F/F^+.$
Recall that $G_{\ell,R}=\prod_{i=1}^{r(R)} GL(m_i(R))_F$ with $\sum_i m_i(R) =m(R).$ The action of the classical Hecke operator $T_{i,j,v}$ for $i = 1, \dots, r(R)$, $j = 1, \dots, m_i(R)$ at unramified place $v$ is recovered by $T_{S,K(R),\CW}(t_{i,j,v})$ with
\[
t_{i,j,v} = \textrm{diag}(\underbrace{\varpi_v,\ldots,\varpi_v}_{j \text{ times}},1,\ldots,1)\times \prod_{k \neq i} \textrm{Id}_{m_k(R)} \in G_{\ell,R,v},
\]
where $\varpi_v$ is the uniformizer and $\textrm{Id}_k$ is the identity $k
\times k$ matrix. The following is a consequence of the main theorem of \cite{Sch}.  

\begin{thm}\label{GaloisrepsR}  Let $\kappa$ be an algebraically closed field and let $\nu:  T_{S,K(R),\CW} \ra \kappa$ be a continuous homomorphism.  Then there are semisimple representations
$$\rho_{i,\nu}:  \Gamma_{F} \ra GL(m_i(R),\kappa), i = 1, \dots, r(R)$$ 
that are characterized, up to equivalence, by the following property:
%\begin{itemize}
%\item[(a)]  
If $v$ is a prime of $F$ not in $S  \cup \{p\}$, then $\rho_{i,\nu}$ is unramified at $v$.  Let $\Gamma_v \subset \Gamma_F$ be a decomposition group at $v$; then the semisimplification of the restriction $\rho_{i,\nu,v}$ of $\rho_{i,\nu}$ to $\Gamma_v$ corresponds to the restriction to $\nu$ to the image of the Hecke operators at $v$ by the unramified Langlands correspondence, in the following sense:  there is an equality of polynomials in $\kappa[X]$
\begin{equation}\label{charpoly}
\det(1 - \rho_{\nu}(\Frob_v)X) =   1 + \sum_{j = 1}^{m_i(R)} (-1)^i \nu(q^{\frac {(m_i(R)+1)j}{2}}T_{i,j,v})X^j
\end{equation}
where the $T_{i,j,v}$ are the standard Hecke operators at $v$ for $GL(m_i)$, normalized as in \cite{Sch}. %?? 
\end{thm}

%{\color{blue} Let $P \subset G$ be a proper rational parabolic subgroup, with fixed Levi subgroup $L_P$.  Then there is an integer $r$, $0 < r \leq n$, and a partition $r = \sum_j r_j$ of $r$, with all $r_j > 0$, such that  
%$$L_P \isoarrow G(V_{n-2r}) \times \prod_{j} GL(r_j)_F,$$
%where $V_{n-2r}$ is a hermitian space over $F$ of dimension $n-2r$ and $G(V_{n-2r})$ is the unitary rational similitude group, defined as in \S \ref{similitude}.   Let $X(L_P)$ denote the locally symmetric space attached to $L_P$, defined as the projective limit over open compact subgroups $U \subset L_P(\af)$ of 
%$$X_U(L_P) = L_P(\QQ)\backslash L_P(\ad)/U \times U_{\infty}Z_{\infty}$$
%for a choice of maximal compact connected subgroup $U_{\infty} \subset L_P(\RR)$, with $Z_{\infty}$ the center of $L_P(\RR)$.
%We let $G_{\ell,P} = \prod_j GL(r_j)_F$. $G_{h,P} = G(V_{n-2r})$.  }  
%{\color{red}  Define local coefficient system on $X_U(L_P)$.}  

\subsection{Eisenstein classes}   The boundary cohomology 
$$\varinjlim_{K} H^*(R\Gamma^{\partial}(S_K(G,X),\CF_K))$$
is expressed in \cite{HZ01} by means of a spectral sequence -- the {\it nerve spectral sequence}
-- whose $E_1$ terms correspond to contributions of the boundary strata corresponding to various standard rational parabolic
subgroups.  The result in \cite{HZ01} is as follows:

\begin{prop}{\cite[Corollary 3.2.9]{HZ01}}\label{nss0}  
There is a spectral sequence abutting to the boundary cohomology in characteristic zero of 
the automorphic vector bundle $\CF$:
\begin{equation} E_1^{r,s}  
%\bigoplus_{r(R) = r+1} \varinjlim_{K_f,\Sigma}
%H^s(\overline{Sh}_{\Sigma}^{R(*)},i^*_R([\mathcal{W}]^{can}))
\Rightarrow \varinjlim_{K} H^{r+s}(R\Gamma^{\partial}(S_K(G,X),\CF_K)).
\end{equation}
The $E_1$ term has a
natural decomposition 
$E_1^{r,s} = \bigoplus_{r(R) = r+1} E_1^{r,s}(R)$,
where $R$ runs over standard rational parabolic subgroups of $G$, and
$$
E_1^{r,s}(R) = Ind_{R(\af)}^{G(\af)}\bigoplus_i \bigoplus_{w \in
W^{R}}E_1^{r,s}(R)_{i,w},$$
where
$$E_1^{r,s}(R)_{i,w} = \widetilde{H}^{s-i-\ell(w)}([\mathcal{F}_{\lambda(h,w)}])\otimes
H^i(Y_\ell(R),\widetilde{\mathbf F}_{\lambda(\ell,w)}).\footnote{The operation $I^R$ that appears in \cite{HZ01} is the identity, because the groups $\Delta_{0,R}$ and $\Delta_{1,R}$ are
trivial in this case, for the reasons explained in \ref{param_boundary}.}
$$
\end{prop}
Here $Y_\ell(R) = \varprojlim_{K(R)} {}_{K(R)}Y_\ell(R)$  is the adelic locally symmetric space attached to the group $G_{\ell,R}$, and  
$\widetilde{H}^{\bullet}([\mathcal{F}_{\lambda(h,w)}])$ denotes coherent cohomology of the Shimura variety attached to $(G_{h,R},X(R))$
with coefficients in an automorphic vector bundle $[\mathcal{F}_{\lambda(h,w)}]$.

Unfortunately, this spectral sequence is not adequate to study the cohomology with $\ZZ_p$ coefficients
at fixed level $K$.  First of all, the decomposition over Weyl group elements $W^{R}$ is based on Kostant's formula for
Lie algebra cohomology, which in general fails in mixed characteristic.   Moreover,  torsion classes at level $K$ cannot in general be identified with the
$K$-invariants in the $G(\af)$-representations $\varinjlim_K H^*(\mathbb{S},\CF_K^{can})$.  So we have to replace
Proposition \ref{nss0} with a calculation of the boundary cohomology at level $K$, and we have to find a less precise 
expression for the individual terms $E_1^{r,s}(R)$ that does not depend on Kostant's formula.  
This is not more difficult but it does require additional notation.   Our approach is roughly analogous to the study of boundary cohomology in \cite[\S\S 3,4]{NT}, with the additional complication that we are working with coherent rather than topological
cohomology, and thus the Levi factor of each parabolic $R$ breaks up into hermitian and linear parts that behave differently.

Denote by
$$i_R:  \partial^R\bbS_{K,\Sigma} \ra \bbS_{K,\Sigma}$$
the closed immersion of the $R$-stratum.   Then we have the {\it nerve spectral sequence}

\begin{equation}\label{nssK} E_1^{r,s} =
\bigoplus_{r(R) = r+1} H^s(\partial^R\bbS_{K,\Sigma},i^*_R(\CF_K^{can}))
\Rightarrow \varinjlim_{K} H^{r+s}(R\Gamma^{\partial}(\bbS_{K,\Sigma},\CF_K))
\end{equation}

Each term in \eqref{nssK} is a module over the Hecke algebra $T_{S,K}$, and the differentials in the spectral sequence are all morphisms of $T_{S,K}$-modules.  We can also let the open compact subgroup $K^p \subset G(\afp)$ shrink to the identity,
letting $K = K_p\times K^p$ with $K_p$ fixed hyperspecial maximal compact, and take the limit over the sequences
\eqref{nssK}, while letting $\Sigma$ vary with $K$ as necessary.   Define
$$R\Gamma^\partial_\infty(\CF_K) = \varinjlim_{K^p} R\Gamma^{\partial}(\bbS_{K,\Sigma},\CF_K)$$
and the $R$-component as
$$R\Gamma^{R}_\infty(\CF_K) = \varinjlim_{K^p,\Sigma} R\Gamma(\partial^R\bbS_{K,\Sigma},i^*_R(\CF_K^{can}))$$
We then obtain a spectral sequence of $\CO[G(\afp)]$-modules
at the limit:

%we let 
%$$\C(R,K) =  R(\QQ)\backslash \pi_0(X) \times G(\af)/K.$$
%We can write 
%$$G(\af) = \coprod_{\alpha \in \cD(R,K)} R(\af)\cdot \alpha \cdot K.$$
%Then the $R$-stratum $\partial_R(\bbS_{K,\Sigma})$ of the toroidal boundary of $\bbS_{K,\Sigma}$
%is the disjoint union of components $\partial_R(\bbS_{K,\Sigma})_\alpha$

\begin{equation}\label{nssin} E_1^{r,s} =
\bigoplus_{r(R) = r+1} H^s(R\Gamma^{R}_\infty(\CF_K))
\Rightarrow  H^{r+s}(R\Gamma^{\partial}_\infty(\CF_K)). 
\end{equation} 

For any finite set $T \supset S$ we let $T_{T,K}$ denote the subalgebra of $T_{S,K}$ generated by the Hecke operators
at primes outside $T$.  We let $G(\ad_f^T) \subset G(\afp)$ denote the subgroup of elements with trivial entries at primes in $T$,
and let $K(T)  = K \cap G(\ad_f^T)$.  The proof of part (iii) of Theorem \ref{Galoisreps}  then comes down to verifying the following three
claims:

\begin{claim}\label{Hecke1}  Let $S'$ be a finite set
of primes containing $S$ and let $\nu:  T_{S',K} \ra k$ be a character realized on a subquotient of $H^{s}(R\Gamma^{\partial}_\infty(\CF_K))\otimes_{\CO} k$ for some $s \geq 0$.  Then there is a proper standard rational parabolic subgroup $R \subset G$, a
finite $T \supset S'$, an $s' \leq s$,
and a character $\nu':  T_{T,K} \ra k$
that coincides with the restriction of $\nu$ and is realized on a subquotient of the space of $K(T)$-invariants of 
$$\varinjlim_{K^p,\Sigma} H^{s'}(\partial^R\bbS_{K,\Sigma},i^*_R(\CF_K^{can}))\otimes_{\CO} k.$$
\end{claim}

In what follows $T' \supset S$ is a finite set as in Claim \ref{Hecke1}.  Assume $K(R,T') = K(T')\cap L_R(\ad_f^{T'})$ is hyperspecial maximal compact at all primes not in $T'$.   We let $T(R)_{T',K}$ denote the product of the unramified Hecke algebras of $L_R(\QQ_v)$
relative to $K(R,T')$ for $v$ not in $T'$ and split in $F/F^+$, and let
\begin{equation}\label{partialsat} s_R:  T_{T',K} \ra T(R)_{T',K}
\end{equation}
be the (unnormalized) partial Satake transform \cite[\S 2.2.6]{NT}.  Here ``unnormalized" means the result of integration along the
unipotent radical of $R$, but without multiplication by the modulus factor $\delta_R^{-1/2}$.

%{\color{red}  Introduce notation for cohomology with local coefficients.}

\begin{claim}\label{Hecke2}  Let $\nu'$ be a character of $T_{T,K}$ as in Claim \ref{Hecke1}.  Then there is an 
automorphic vector
bundle $\CF^R$ on $\bbS_{K_R}(G_{h,R},X(R))$, a finite-dimensional representation $\CW$ of $G_{\ell,R}$ with coefficients
in a finite extension $k'/k$, a finite set $T' \supset T$
as above such that $K(T')\cap L_R(\ad_f^{T'})$ is hyperspecial maximal at all primes not in $T'$, and a
character $\nu'_R:  T(R)_{T',K} \ra k'$ such that
\begin{itemize}
\item[(a)]  $\nu'_R$ is realized on a subquotient of 
$$H^i(\bbS_{K_R}(G_{h,R},X(R))_{\Sigma_R},\CF^{R, \rm{can}})_k\otimes H^j({}_{K(R)}Y_\ell(R),\tilde{\CW})$$
for some $i, j \geq 0$;
\item[(b)]$\nu'_R \circ s_R$ coincides with the restriction of $\nu'$ to $T_{T',K}$.
\end{itemize}
\end{claim}

%{\color{red}  In \S \ref{Hecke2proof} give more details on (a) and say something about (b)}

\begin{claim}\label{Hecke3}  Notation is as in Claim \ref{Hecke2}.  Let $\nu'$ be a character of $T_{T,K}$ with values 
in a finite extension $k'$ of $k$. Suppose $R$, $T'$, and $\nu'_R$ are as in Claim \ref{Hecke2}.  Then the residual Galois 
representation $\br_{\nu'}$ attached to $\nu'$ by Theorem \ref{Galoisreps} (i) is reducible.
\end{claim}

\subsubsection{Proof of Claim \ref{Hecke3}}

The proofs of Claims \ref{Hecke1} and \ref{Hecke2} are postponed until the next section.  Assuming both of these
Claims,  \ref{Hecke3} can now be proved right away. 
The bulk of the work is already contained in \S 4 of \cite{NT}.  We introduce some additional notation.  
 By Claim \ref{Hecke1} we may assume $K(R,T')$ factors as a product $K_{h,R,T'}\times K_{\ell,R,T'}$
with $K_{h,R,T'} = K(R,T')\cap G_{h,R}(\ad_f^{T'})$, $K_{\ell,R,T'} = K(R,T')\cap G_{\ell,R}(\ad_f^{T'})$.
Let 
$$T(R)_{h,T',K} \text{   (resp. $T(R)_{\ell,T',K}$)}$$
denote the product of the unramified Hecke algebras of $G_{h,R}(\QQ_v)$ (resp. $G_{\ell,R}(\QQ_v)$) relative to $K_{h,R,T'}$ (resp. $K_{\ell,R,T'}$)
for $v$ not in $T'$ and split in $F/F^+$.   The product map $G_{h,R} \times G_{\ell,R} \isoarrow L_R$ defines a canonical isomorphism
\begin{equation}\label{phl}
 \varphi_R:  T(R)_{T',K} \isoarrow T(R)_{h,T',K}\otimes T(R)_{\ell,T',K} .
 \end{equation}
Let 
$T^i_{h,K(R),T',\CF^R,k}$ (resp. $T^j_{\ell,K(R),T',\CW}$) denote the image of the natural map
$$T(R)_{h,T',K} \ra \Endo(H^i(\bbS_{K_R}(G_{h,R},X(R))_{\Sigma_R},\CF^{R, \rm{can}})_k)$$
(resp. 
$$T(R)_{\ell,T',K} \ra \Endo(H^j({}_{K(R)}Y_\ell(R),\tilde{\CW})).$$
Then (possibly replacing $k'$ by a finite extension) there are characters
$$\nu'_h:  T^i_{h,K(R),T',\CF^R} \ra k'; ~~ \nu'_\ell\ra k':  T^j_{\ell,K(R),T',\CW}$$
such that 
\begin{equation}\label{phl2}
\nu' = \nu'_h\otimes \nu'_\ell \circ (\varphi_R\otimes k').
\end{equation}
Now it follows from Theorem \ref{Galoisreps} (i), applied to $\nu'_h$, that there exists an 
$n_h := n- 2m(R)$-dimensional representation $\br(\nu'_h)$ of $\Gamma_F$ such that, for all $v$ split in $F/F^+$ and
outside $T'$, 
\begin{equation}\label{charpolyh}
\det(1 - \br(\nu'_h)(\Frob_v)X) =   1 + \sum_{j = 1}^{n_h} (-1)^j \nu(q^{\frac {(n_h+1)j}{2}}T_{j,v} )X^j
\end{equation}
Similarly, with notation as in Theorem \ref{GaloisrepsR}, 
  there are $m_i(R)$-dimensional representations
$\br_i(\nu'_\ell)$ of $\Gamma_F$ with coefficients in $k'$, $i = 1, \dots, r(R)$, such that, for all $v$ split in $F/F^+$ and
outside $T'$, 
\begin{equation}\label{charpolyell}
\det(1 - \br_i(\nu'_\ell)(\Frob_v)X) =   1 + \sum_{j = 1}^{m_i} (-1)^j \nu(q^{\frac {(m_i+1)j}{2}}T_{i,j,v} )X^j.
\end{equation}

Let $\br_i^d(\nu'_\ell) = \br_i(\nu'_\ell)^{\vee}(1-m_i(R))$ (Tate twist).  It then follows from identity (b) of Claim \ref{Hecke2}, as in the 
discussion following \cite[Lemma 4.6]{NT}, that there are integers $\mu_0$ and $\mu_1^{\pm},\mu_2^{\pm}, \dots, \mu_{r(R)}^{\pm}$
% By Claim \ref{Hecke2} and Chebotarev density (which allows us to increase $T'$ if necessary) it suffices to replace $\nu'$ by $\nu'_R \circ s_R$ in the statement of Claim \ref{Hecke3}.
such that $\nu'_R$ and 
$$\br(\nu'_h)(\mu_0) \oplus \sum_{i = 1}^{r(R)} \br_i(\nu'_R)(\mu_i^+) \oplus \br_i^d(\nu'_R)(\mu_i^-),$$
(where $(\mu_i^{\pm})$ denote Tate twist)
have the same  Frobenius eigenvalues outside of $T'$. More precisely, these are given by $\mu_0 = n_h - n = -2m(R)$ and
\begin{equation}\label{Tatetwist}
        \mu_i^+ = -m_i(R) - 2\sum_{j>i} m_j(R),\quad \quad 
        \mu_i^- = 2\sum_{j>i} m_j(R)
\end{equation}
for $i = 1,\ldots, r(R).$
By Chebotarev density, this completes the proof.

%\subsubsection{Coefficients on boundary spaces}  

\subsubsection{Eisenstein characters}

Let $k$ be an algebraically closed field (of any characteristic).  Let $R \subset G$ be a rational parabolic subgroup.
Let $\CS$ be a finite set of primes of $F^+$, and  let $\CT_R^\CS$ be the restricted tensor product of the local Hecke algebras of
$L_R$, with coefficients in $k$, at all primes of $F^+$ outside $\CS$ that split in $F$; let $\CT^\CS = \CT_G^\CS$.  Here by
{\it restricted} we mean that the local component is the identity at all but finitely many places.
Choose a maximal split torus $T_v \subset L_R$ for each prime $v$ of $F^+$ that splits in $F$, and let
$\CT_{T_v}$ denote its Hecke algebra over $k$.  The Satake homomorphism is an injective map
\begin{equation}\label{satG}
s_G:  \CT^{\CS}  \ra \otimes'_{v \notin \CS} \CT_{T_v}
\end{equation}
(again, $v$ runs only over primes that split in $F$).
Similarly, we have a Satake homomorphism for $L_R$
\begin{equation}\label{satR}
s_R:  \CT_R^{\CS}  \ra \otimes'_{v \notin \CS} \CT_{T_v}
\end{equation}

Corresponding to the factorization $L_R \isoarrow G_{h,R} \times G_{\ell,R}$,
we write $\CT_R^{\CS} = \CT_{h,R}^{\CS} \otimes \CT_{\ell,R}^{\CS}$.   
The character $\beta:  \CT_R^{\CS} \ra k$ is called {\it cohomological} if it is the tensor product of characters
$\beta_h$ of $\CT_{h,R}^{\CS}$ and $\beta_\ell$ of $\CT_{\ell,R}^{\CS}$, where $\beta_h$ occurs as a subquotient of 
the coherent cohomology of a (smooth projective) toroidal compactification of the integral model $\bbS_{K_R}(G_{h,R},X(R))$ with coefficients in some canonically extended automorphic vector bundle, for some level subgroup $K_R$ that is hyperspecial
maximal compact at split primes not in $\CS$, and where $\beta_\ell$ occurs as a subquotient of the cohomology (with some 
coefficients) of the locally symmetric space attached to $G_{\ell,R}$, with full level at split primes not in $\CS$.

\begin{defn}\label{eis}  Let $\nu:  \CT^{\CS} \ra k$ be a character of $\CT^{\CS}$.  We say $\nu$ is {\it Eisenstein}
if, viewing it via \eqref{satG} as a character of $s_G(\CT^{\CS}) \subset \otimes'_{v \notin \CS} \CT_{T_v}$, it extends
to a character $\beta': s_R(\CT_R^{\CS}) \ra k$ so that $\beta = \beta' \circ s_R$ is cohomological in the sense 
just defined.
\end{defn} 

Let $K_f \subset G(\af)$ be a (neat) level subgroup, and let $K_R \subset K_f \cap L_R(\af)$ be a compact open subgroup.  

\subsection{Proof of Claim \ref{Hecke1}}

Let $K' \subset K \subset G(\af)$, with $K'_p = K_p$ our fixed hyperspecial maximal compact and $K'$ normal in $K$.  Then for any fixed toroidal datum $\Sigma$, we have a finite morphism
$$\bbS_{K',\Sigma} \ra \bbS_{K,\Sigma}$$
of normal schemes.  We assume $\bbS_{K,\Sigma}$ to be smooth and projective for convenience.  The proof of \cite[Lemma 2.6]{H89} applies in the mixed characteristic situation and implies that
\begin{lem}  The above map defines a canonical isomorphism 
$$\bbS_{K',\Sigma}/(K/K') \isoarrow \bbS_{K,\Sigma}.$$
In particular, for any automorphic vector bundle $\CF_{\bullet}$ there are canonical isomorphisms (Hochschild-Serre spectral sequence) in the derived category of $\CO$-modules
$$RHom_{K/K'}(\CO,R\Gamma(\bbS_{K',\Sigma},\CF_{K'}^{\rm can})) \isoarrow R\Gamma(\bbS_{K,\Sigma},\CF_{K}^{\rm can});$$
$$RHom_{K/K'}(\CO,R\Gamma(\bbS_{K',\Sigma},\CF_{K'}^{\rm sub})) \isoarrow R\Gamma(\bbS_{K,\Sigma},\CF_{K}^{\rm sub});$$
$$RHom_{K/K'}(\CO,R\Gamma^{\partial}(\bbS_{K',\Sigma},\CF_{K'}^{\rm can})) \isoarrow R\Gamma^{\partial}(\bbS_{K,\Sigma},\CF_{K}^{\rm can}).$$
\end{lem}

\begin{proof}  The first claim, as noted, follows as in the proof of  \cite[Lemma 2.6]{H89}.  The first two spectral sequences
then follow from Proposition \ref{db} (iii), and the third from the definition of $R\Gamma^{\partial}$.
\end{proof}

\subsubsection{Proof of Claim \ref{Hecke1}}
Passing to the limit over $K'$ (containing $K_p$), the last isomorphism yields a Hochschild-Serre spectral sequence of $\CO$-modules:
\begin{equation}\label{hsss}
E_2^{a,b} = H^a(K^p,H^b(R\Gamma^\partial_\infty(\CF_K))) \Rightarrow H^{a+b}(R\Gamma^{\partial}(\bbS_{K,\Sigma},\CF_{K}^{\rm can})).
\end{equation}
For $T$ as in the statement of Claim \ref{Hecke1}, let $K^p_T$ denote the subgroup of elements of $K^p$ whose entries at each prime $v$ in $T$ belongs to the principal congruence subgroup $K(v) \subset K_v$ of elements congruent to the identity modulo $v$,
$K_T = K^p_T \times K_p$.
Since $K(v)$ has pro-order prime to $p$, we can rewrite 
\begin{equation}\label{hsssK} H^a(K^p,H^b(R\Gamma^\partial_\infty(\CF_K))) = 
\varinjlim_{T \supset S'} H^a(K^p/K^p_T, H^b(R\Gamma^\partial_\infty(\CF_K))^{K_T})
\end{equation}
(compare \cite[Lemma 4.2]{NT}). For each $T$, the Hecke algebra $T_{T,K_T}$ acts on the $K_T$-invariants on the left-hand side compatibly with the action on 
$H^\bullet(R\Gamma^\partial_\infty(\CF_K))^{K_T}$ on the right hand side.  Claim \ref{Hecke1} now follows by combining
\eqref{hsssK} with \eqref{nssin}.

%\subsection{Nerve spectral sequence} 

%\begin{equation}\label{nss} E_1^{r,s} =
%\bigoplus_{r(R) = r+1} \varinjlim_{K_f,\Sigma}
%H^s(\overline{Sh}_{\Sigma}^{R(*)},i^*_R([\mathcal{W}]^{can}))
%\Rightarrow \varinjlim_{K} H^{r+s}(R\Gamma^{\partial}(S_K(G,X),\CF_K))
%\end{equation}
%Here $H^{\bullet}([\mathcal{W}](\infty))$ denotes
%$\varinjlim_{K_f,\Sigma} H^{\bullet}(\partial Sh_{\Sigma},
%[\mathcal{W}]^{can}\otimes \mathcal{O}_{\partial Sh_{\Sigma}})$,

\section{Automorphic vector bundles on the toroidal boundary}\label{boundary}

The calculation of the boundary coherent cohomology in \cite{HZ94,HZ01} uses a combination of algebraic,
analytic, and representation-theoretic arguments that do not apply to the integral models studied in this paper.
In fact, some of the arguments probably fail in small characteristic.  For example, if $f:  A \ra S$ is the canonical morphism of a
universal abelian scheme $A$ (with some level structure) to a Shimura variety $S$ of PEL type, the computation in \cite{HZ94} of the 
higher direct images $R^if_*$ of an automorphic vector bundle on $A$ as automorphic vector bundles on $S$ is based
on Kostant's theorem on Lie algebra cohomology of unipotent radicals of parabolic subalgebras of a reductive Lie algebra,
and this just breaks down in general.  

Fortunately, most of what we need has been worked out by Kai-Wen Lan for the toroidal compactifications of \cite{Lan}.  As in  \cite{HZ94,HZ01}, the calculation is carried out for individual closed strata of the boundary, which gives coherent cohomology of
the boundary strata attached to a Shimura datum $(G_{h,R},X(R))$ in the minimal compactification; these are then put together according to the configuration of the strata, which introduces the topological cohomology of the factor $G_{\ell,R}$.
The first important observation is that canonical extensions of automorphic vector bundles behave well under restriction to
boundary strata:  this corresponds to the discussion around Corollary 3.4.3 of \cite{HZ94} and the first part of Proposition 5.6
of \cite{Lan4}.

\subsection{The algebraic part}\label{algebraicpart}

To begin, we write  $\bbS_{K,\Sigma}$ for the integral model of the toroidal compactification, previously
denoted $\bbS_K(G,X)_\Sigma$.    We always choose $\Sigma$ as in Theorem \ref{toroidalintegral},
so that the compactification is proper and smooth and the boundary divisor is a divisor with smooth normal crossings.
Fix a rational boundary component $F$ as in \S \ref{termin} -- we may assume it corresponds to (a component of)
a Shimura datum $(G_{h,R},X(R))$ with $R=P_F$, its stabilizer; in other words $F \subset X(R)$ -- let $\sigma \in \Sigma_F$ be a cone and let 
$$i_\sigma: Z_\sigma \hookrightarrow \bbS_{K,\Sigma}$$
 denote the inclusion of the corresponding {\it locally closed} stratum of the toroidal boundary,  $\bar{Z}_\sigma$ its closure in $\bbS_K(G,X)_\Sigma$.  Then we have a diagram as in \cite[(1.2.5)]{HZ94}, \cite[(4.4),(4.5)]{Lan4} of smooth
 schemes over $Spec(\CO)$:
\begin{equation}\label{2step}
Z_\sigma \overset{\pi_2}\to  A_F \overset{\pi_1} \to M_F \subset \bbS_{K'}(G_{h,R},X(R))
\end{equation}
for some appropriate level subgroup $K' \subset G_{h,R}(\af)$, where $\pi_2$ is a torus fibration and $\pi_1$ is an abelian
scheme over a connected component $M_F$ of $\bbS_{K'}(G_{h,R},X(R))$.

Now let $\CF_K^{\rm can}$ be a canonically extended automorphic vector bundle over $\bbS_{K,\Sigma}$.

\begin{prop}\cite[Proposition 5.6]{Lan4}\label{pullback1}  (i)  There is a canonically determined locally free coherent
sheaf $\CF^A_K$ over $A_F$
and a canonical isomorphism
$$i_\sigma^*(\CF_K^{\rm can}) \isoarrow \pi_2^*(\CF^A_K)$$
compatible with the action of $G(\afp)$ (see below for explanation).

(ii)  The vector bundle $\CF^A_K$ is endowed with an increasing  $P_F(\afp)$-invariant  filtration by vector bundles
$$\dots \subset \CF^{A,j}_{K} \subset \CF^{A,j+1}_{K} \subset \dots$$
such that each associated graded piece $\gr^j(\CF^{A,j}_{K})$ is isomorphic to the pullback via $\pi_1$
of an automorphic vector bundle on $\bbS_{K'}(G_{h,R},X(R))$:
$$\gr^j(\CF^{A,j}_{K}) \isoarrow \pi_1^*(\CF_{K',R}^j).$$
\end{prop}

\begin{rmk} (i) The group $G(\afp)$ acts on the set of toroidal compactifications $\bbS_{K,\Sigma}$ by acting on the
group $K$ as well as the toroidal datum $\Sigma$.  The subgroup $P_F(\afp) \subset G(\afp)$ fixes the set of toroidal boundary components, as $K$ and $\Sigma$ vary, over the rational boundary component $F$, and the isomorphisms are compatible with diagrams \eqref{2step}
and thus with the morphisms $\pi_1, \pi_2$ (which can be labelled with $\sigma$).  Detail are left to the reader, but
see the discussion on pp. 310-11 of \cite{HZ94} and p. 89 of \cite{HZ01}.  

(ii)  Part (ii) of Proposition \ref{pullback1} corresponds to the discussion on pp. 325-6 of \cite{HZ94}.  In \cite{Lan4} it
is not stated explicitly that the graded pieces correspond to automorphic vector bundles on the base $M_F$,
but the proof makes it clear that the ``coherent sheaves over $\mathbf{Z}$" of the statement of Lan's Proposition 5.6
are indeed the automorphic vector bundles as constructed by Lan in \cite{Lan2}.

The $P_F(\afp)$-invariance of the filtration follows from the fact that the construction of automorphic vector bundles in \cite[Definition 6.7]{Lan2} is a functor from locally free $\CO$-representations of $P_F$ to $P_F(\afp)$-equivariant vector bundles over $A_F$.  The $P_F(\afp)$-action is not mentioned explicitly in \cite{Lan2} but the construction is exactly analogous to that in \cite[Proposition 8.1.4.1]{Lan3}.
\end{rmk}

\begin{prop}\label{lanextension}  After replacing $\Sigma$ by a refinement if necessary, the morphism $\pi_2$ of \eqref{2step} extends to a morphism 
$$\bar{\pi}_2 = \bar{\pi}_{2,\sigma}:  \bar{Z}_\sigma \ra A_{F,\Sigma'}$$
where $A_{F,\Sigma'}$ is an appropriate (smooth projective) toroidal compactification of $A_F$.
Moreover, letting $\bar{i}_\sigma:  \bar{Z}_\sigma \hookrightarrow \bbS_{K,\Sigma}$ denote the inclusion of the closed
boundary stratum, the isomorphism of Proposition \ref{pullback1} extends to a canonical isomorphism
$$\bar{i}_\sigma^*(\CF_K^{\rm can}) \isoarrow \bar{\pi}_2^*(\CF^{A,\rm can}_K).$$
\end{prop}

\begin{proof}  The first part is implicit in \cite{Lan3}; a complete proof  will appear in future work of Lan.  To prove the second part, note that
$\bar{i}_\sigma^*(\CF_K^{\rm can})$ and $\bar{\pi}_2^*(\CF^{A,\rm can}_K)$ are both vector bundles on the regular scheme
$\bar{Z}_\sigma$ that are known to be canonically isomorphic away from a subscheme of codimension $2$ by Proposition \ref{pullback1} and by Proposition 3.12.2 of \cite{HZ94}, together with the results of \cite{Lan12}.   Thus there 
is a canonical isomorphism over $\bar{Z}_\sigma$.  
\end{proof}

\begin{prop}\cite{Lan2}\label{pi1extension} %{\color{red}  Lan reference?} 
Let $A_F \hookrightarrow A_{F,\Sigma'}$ be the smooth projective toroidal compactification of $A_F$
of Proposition \ref{lanextension}.  After possibly replacing $\Sigma'$ by a refinement, say $\Sigma''$, the map $\pi_1$ of
\eqref{2step} extends to a morphism
$$\bar{\pi}_1:  A_{F,\Sigma''} \ra \bar{M}_F \subset  \bbS_{K'}(G_{h,R},X(R))_{\Sigma_F}$$
where $\bbS_{K'}(G_{h,R},X(R))_{\Sigma_F}$ is a smooth projective toroidal compactification of $\bbS_{K'}(G_{h,R},X(R))$
and $\bar{M}_F$ is the Zariski closure of $M_F$.
\end{prop}
\begin{proof}  The main reference for this fact is Theorem 2.15 of \cite{Lan2}.  Since the notation $A_{F,\Sigma'}$ is not 
quite justified there, Lan advises us to add a reference to Section 1.3.4 of \cite{Lan3}.  Although the title of \cite[\S 1.3]{Lan3} is 
``Algebraic compactifications in characteristic zero," the constructions in the relevant section are valid in mixed characteristic.
\end{proof}

We use the same notation $\CF^{A,\rm can}_K$ to denote the canonical extension of $\CF^A_K$ on $A_{F,\Sigma''}$.  

\begin{prop}\cite{Lan2}\label{directimages}  (i)  The filtration $\{\CF^{A,j}_{K} \}$ of $\CF^A_K$ of Proposition \ref{pullback1} (ii) 
extends to a filtration $\{\CF^{A,j}_{K,\Sigma'} \}$ of $\CF^{A,\rm can}_K$ by vector bundles such that each
$\gr^j(\CF^{A,\rm can}_K)$ is the pullback via $\bar{\pi}_1$ of the canonical extension of  $\CF_{K',R}^j$.

(ii)  The higher direct images $R^k{\bar{\pi}_1}\gr^j(\CF^{A,\rm can}_K)$, $k \geq 0$, are canonical extensions of the automorphic vector bundles
$R^k\pi_1\gr^j(\CF^{A}_K)$ on $M_F$.
\end{prop}

\begin{proof}   (i)  We thank Lan for the following argument.  The point is that the construction in \S 3.B of \cite{Lan2} defines
the canonical extensions algebraically in terms of the relative differentials on semiabelian schemes over $M_F$.   
%{\color{blue}  Add details to the following.}  
More precisely, the
semi-abelian scheme over the toroidal compactification
$A_{F, \Sigma'}$, denoted $\tilde{G}$ in {\it loc. cit.},
admits a  
split subtorus $\CT$ and hence a quotient
semi-abelian scheme $\bar{G}$.   It is explained in {\it loc. cit.} that this $\bar{G}$ is the
pullback via $\bar{\pi}_1$ from a semi-abelian scheme (denoted $G$ in {\it loc. cit.} over $\bar{M}_F$ .  In
particular, since Lan shows that the canonical extension of a given automorphic vector bundle on $A_{F, \Sigma'}$ arises
as a subquotient of some tensor power of the relative
Lie algebra of $\tilde{G}$,  it admits an extension structure coming
from the short exact sequence
$$1 \ra \CT \ra \tilde{G} \ra \bar{G} \ra 1$$
of group schemes over $A_{F,\Sigma'}$.  
This structure induces filtrations on all the canonical extensions in such as way as to guarantee that the graded pieces are defined 
by the relative Lie algebras of $\CT$ and $\bar{G}$; but each of these is a pullback from the toroidal compactification $\bar{M}_F$.

%{\color{red}  (i) Find references.}

(ii)  We first observe that each $R^k{\bar{\pi}_1}\gr^j(\CF^{A,\rm can}_K)$ is a locally free sheaf on $\bar{M}_F$.   Indeed, by (i),
each $\gr^j(\CF^{A,\rm can}_K)$ is a pullback via $\bar{\pi}_1$, so by the projection formula the claim reduces to the corresponding
assertion for the structure sheaf of $A_{F,\Sigma''}$, where it follows from \cite[Theorem 2.15]{Lan2}.   The same result of
Lan implies that our statement is true outside a divisor on the special fiber, and we know by the results of \cite{HZ94} that the statement is true on the generic fiber, and is true outside
a divisor on the special fiber by results of Lan.  Thus, 
as in the proof of Proposition \ref{lanextension}, it is true globally, since $\bar{M}_F$ is a regular scheme.
\end{proof}

%%%%%%%%%%%%%%%%%%%%%%%%%%%%%%%%%%%%%%%%%%%%%%%%
%%%%%%%%%%%%%%%%%%%%%%%%%%%%%%%%%%%%%%%%%%%%%%%%%
\subsection{The topological part}\label{topologicalpart}

In order to derive Claim \ref{Hecke2} from Claim \ref{Hecke1} we need to understand 
$H^{s'}(\partial^R\bbS_{K,\Sigma},i^*_R(\CF_K^{can}))\otimes_{\CO} k$ as $K^p$ and $\Sigma$ vary.
We need to fix $\Sigma$ (and therefore $K^p$) in order to state the next Proposition.

\begin{prop}\label{311HZ3}

\begin{itemize}
\item[(i)] For each $b \geqslant 0$, the assignment
$$
\sigma \longmapsto H^{b}\left(\bar{Z}_{\sigma}, i_{\sigma}^{*}(\CF_K)^{can}\right)
$$
defines a locally constant sheaf $\boldsymbol{L}^{b}(\cdot, \CF_K)$ of $\CO$-modules 
on the simplicial complex $\mathfrak{N}_{\Sigma}(R)$, and an associated spectral sequence
$$
 \quad E_{1}^{a, b}=H^{a}\left(\mathfrak{N}_{\Sigma}(R), \boldsymbol{L}^{b}(\cdot, \CF_K)\right) \Longrightarrow H^{a+b}\left(\partial^R\bbS_{K,\Sigma}, i_{R}^{*}(\CF_K)^{can}\right).
$$
\item[(ii)]  Similarly, for each $b \geqslant 0$, the assignment

$$
\sigma \longmapsto H^{b}\left(\bar{Z}_{\sigma}, \bar{\pi}_{2,\sigma}^*(\CF_{K,\Sigma'}^{A,j})\right)
$$ 
defines a Hecke-equivariant locally constant sheaf 
$\boldsymbol{L}^{b}(\cdot, \CF^j_{K'})$ of $\CO$-modules, written $\sigma \mapsto  \boldsymbol{L}^{b}(\cdot, \CF^j_{K'})(\sigma)$ with $K'$ as in Propositions~\ref{pullback1},\ref{pi1extension},
on the simplicial complex $\mathfrak{N}_{\Sigma}(R)$.
There is a Hecke-equivariant spectral sequence derived from the filtration in Proposition \ref{directimages} 
$$E_1^{a',j} = H^{a'}(\mathfrak{N}_{\Sigma}(R),\boldsymbol{L}^{b}(\cdot, \CF^j_{K'})) \Rightarrow H^{a'+b+j}\left(\mathfrak{N}_{\Sigma}(R), \boldsymbol{L}^{b}(\cdot, \CF_K)\right). $$

\item[(iii)]  For each $\sigma$ there is a canonical Leray spectral sequence
$$E_2^{s,t} = H^s(\overline{M}_F,R^t\bar{\pi}_{1,*}\CO_{A_{F,\Sigma''}}\otimes \CF^j_{K'}) \Rightarrow
\boldsymbol{L}^{s+t}(\cdot, \CF^j_{K'})(\sigma)
$$
\item[(iv)]  Suppose $\sigma$ is a face of $\sigma'$, and let
$$i_{\sigma,\sigma'}:  \bar{Z}_{\sigma'} \hookrightarrow \bar{Z}_\sigma$$
denote the corresponding closed immersion.  Then for each $j$ the pullback map
$$i_{\sigma,\sigma'}^*:   \bar{\pi}_{2,\sigma'}^*(\CF_{K,\Sigma'}^{A,j}) \ra  \bar{\pi}_{2,\sigma}^*(\CF_{K,\Sigma'}^{A,j})$$
determines isomorphisms on cohomology
$$\boldsymbol{L}^b(\cdot,\CF^j_{K'}))(\sigma) \ra \boldsymbol{L}^b(\cdot,\CF^j_{K'}))(\sigma')$$
that are the face maps for the local system on  $\mathfrak{N}_{\Sigma}(R)$.

%{\color{red}  EXPLAIN the maps as $\sigma$ varies}
%, and an associated spectral sequence
%$$
% \quad E_{1}^{a, b}=H^{a}\left(\mathfrak{N}_{\Sigma}(R), %\boldsymbol{L}^{b}(\cdot, \CF^j_{K'})\right) \Longrightarrow %H^{a+b}\left(\partial^R\bbS_{K,\Sigma}, i_{R}^{*}(\CF_K)^{can}\right).
%$$
\end{itemize}

\end{prop}
\begin{proof} Part (i) is proved as in \cite{HZ94,HZ01}.  We compute the cohomology by the spectral sequence for the closed cover of $Z_\Sigma(R)$
by  the maximal strata $\bar{Z}_\sigma$ as $\sigma$ runs over the simplices in $\mathfrak{N}_{\Sigma}(R)$.
Proposition \ref{lanextension} allows us to apply the proof of \cite[Proposition 3.1.1 (i)]{HZ01}.   Part (ii)  is a consequence  Proposition \ref{directimages}.
Part (iii) follows from the projection formula, in view of the isomorphism 
\begin{equation}\label{iso-boundary-abelian}
H^{b}\left(\bar{Z}_{\sigma}, i_{\sigma}^{*}(\CF_K)^{can}\right)\isoarrow H^b(A_{F,\Sigma''}, \CF_K^{A,can}),
\end{equation}
which is proved as \cite[Proposition 3.4.1]{HZ01} in light of the results of this section
and
$$H^b(A_{F,\Sigma''},\gr^j\CF^{A,can}_{K}) \isoarrow H^b(A_{F,\Sigma''},\bar{\pi}_1^*(\CF_{K',R}^{j,can})).$$

  Part (iv) is just a restatement of the fact that the functor 
$\sigma \mapsto \boldsymbol{L}^b(\cdot,\CF^j_{K'}))(\sigma)$ defines a locally constant sheaf on $\mathfrak{N}_{\Sigma}(R)$.
\end{proof}

Combining the three parts of Proposition \ref{311HZ3} we find:

\begin{cor}\label{reductionRh}  Let $\nu'$ be a character of $T_{T,K}$ as in Claim \ref{Hecke1}.  In order to prove Claim \ref{Hecke2}, it suffices to prove it if $\nu'$ is realized
in the action of $T_{T,K}$ (after possibly increasing $T$) on $H^{a}(\mathfrak{N}_{\Sigma}(R),Tor_c^\CO(\boldsymbol{L}^{b}(\cdot, \bar{\pi}_1^*\CJ_{K'}),k))$
for an automorphic vector bundle $\CJ_{K'}$ on $M_F$ and for $c = 0, 1$, where the cohomology of $\bar{\pi}_1^*\CJ_{K'}$ is computed via (\ref{iso-boundary-abelian}).  
\end{cor}

%%%%%%%%%%%%%%%%%%%%%%%%%%%%%%%%%%%%%%%%%%%%%%%%
%%%%%%%%%%%%%%%%%%%%%%%%%%%%%%%%%%%%%%%%%%%%%%%%%
\subsubsection{The homotopy type of $\mathfrak{N}_{\Sigma}(R)$}\label{homotopynerve}

\begin{prop}\label{homotopy}  The nerve $\mathfrak{N}_{\Sigma}(R)$ of the $R$-stratum $\partial^R\bbS_{K,\Sigma}$ is homotopy equivalent to  the $G_{\ell,R}$-stratum of the  Borel-Serre compactification of the locally symmetric space ${}_{K(R)}Y_\ell(P(R))$.
\end{prop}
%{\color{red}  Restate with level subgroups}
\begin{proof}  This is \cite[Proposition 2.6.4]{HZ94b}.
\end{proof}

%%%%%%%%%%%%%%%%%%%%%%%%%%%%%%%%%%%%%%%%%%%%%%%%
%%%%%%%%%%%%%%%%%%%%%%%%%%%%%%%%%%%%%%%%%%%%%%%%%
\subsubsection{Local systems on adelic and connected locally symmetric spaces}
We let $H$ be a connected reductive algebraic group over $\QQ$ with center $Z$, let $K_H \subset H(\RR)$ be a maximal connected subgroup that is compact modulo $Z(\RR)$, and let $Y = H(\RR)/K_H$ be the corresponding (possibly disconnected) symmetric space.  Let $K_f \subset H(\af)$ be a compact open subgroup and let
$$S(Y,K_f) = H(\QQ)\backslash Y \times H(\af)/K_f$$
denote the adelic symmetric space.  Write $H(\af) = \coprod_j H(\QQ)\alpha_j K_f$ and let $\G_j = H(\QQ)\cap \alpha_j K_f \alpha_j^{-1}$, $S(\G_j) = \G_j\backslash Y$.
Then we can rewrite
\begin{equation}\label{componentsK} S(Y,K_f) = \coprod_j S(\G_j).
\end{equation}

Now let $M$ be a finite abelian group, and let $\lambda:  K_f \ra \Aut(M)$ be a homomorphism.  Consider the local system over $S(Y,K_f)$
$$M(K_f) = H(\QQ)\backslash Y \times H(\af)\times M/K_f$$
where $K_f$  (resp. $H(\QQ)$) acts diagonally on the last two (resp. first two)  factors in the product $Y \times H(\af)\times M$.  

\begin{prop}\label{localsysadelic}  Under these hypothesis, the restriction $M(K_f)_j$ of the local system $M(K_f)$ to the subspace $S(\G_j) \subset S(Y,K_f)$  can be written as
\[
\G_j\backslash Y\times M,
\]
where $\G_j$ acts via $q\cdot (y,m) = (qy, \lambda(\alpha^{-1}_j q^{-1}\alpha_j) \cdot m)$.
\end{prop}
\begin{proof}
 For $\alpha_j$ as above, consider 
 \[
 [y,m]\mapsto [y,\alpha_j,m]: \G_j\backslash Y\times M\to H(\QQ)\backslash Y \times H(\af)\times M/K_f
 \]
 To assert injectivity note that $[y,\alpha_j,m]=[y',\alpha_j,m']$, then $y=qy',~\alpha_j = q\alpha_j k$ and $\lambda(k) m' =m$ for some $q
 \in H(\QQ)$ and $k\in K_f$. From the second equality we obtain that $q\in \G_j$ and $k = \alpha^{-1}_j q^{-1}\alpha_j$, so that $[y,m]=[y',m'].$ 
 
 To show that $M(K_f)$ is the disjoint union of the local systems $M(K_f)_j$, pick $(y,h,m)\in Y \times H(\af)\times M$. Then $h = q\alpha_jk$ for some $q\in H(\QQ)$ and $k\in K_f$. Then $[y,h,m] = [q^{-1}y,\alpha_j, \lambda(k)^{-1} m]$, which is in the image of $M(K_f)_j.$ Lastly, if $[y,\alpha_j, m]=[y',\alpha_i, m'],$ then $y =qy',~\alpha_j = q\alpha_i k,$ and $m = \lambda(k)m'$. The second equation implies $i=j.$
\end{proof}

Let $H = G^{red}_{\ell,R}$, $Y_{\ell,R}$ be the corresponding symmetric space, and $\{\alpha_j\}$ as above. Let $\bold{\lambda}:  H \ra \MM$ be an irreducible algebraic representation defined over $Spec(\Zp)$.   Fix an $H(\Zp)$-invariant lattice
$\Lambda \subset \MM(\Qp)$ and let $M = \Lambda/p\Lambda$.  The $H(\Zp)$-module $M$ may depend on the lattice but by the Brauer-Nesbitt theorem
its semisimplification does not.   Applying Proposition \ref{localsysadelic} to this situation, we find

\begin{cor}\label{lattices}  Let $\tilde{\MM}$ be the local system 
$$\tilde{\MM} = H(\QQ)\backslash Y_{\ell,R} \times H(\af)\times \MM(\Qp)/K_f.$$
Then there is a local system of free $\Zp$-modules $\tilde{\Lambda} \subset \tilde{\MM}$ with the property that, for each component $S(\G_j)$ as above, the restriction
of $\tilde{\Lambda}/p\tilde{\Lambda}$ to $S(\G_j)$ is isomorphic to
\[
\G_j\backslash (Y\times \tilde{\Lambda}/p\tilde{\Lambda}),
\]
with the twisted $\G_j$-action via $q\cdot (y,\tilde{m}) = (qy, \lambda(\alpha^{-1}_j q^{-1}\alpha_j)\cdot \tilde{m})$.
\end{cor}

Consider a compact open $K\subseteq G(\af)$ satisfying the properties listed in Appendix~\ref{param_boundary}. For a standard rational parabolic $R$ with Levi decomposition $R=L_R\cdot U_R$, let $\Gamma_{R} \coloneqq R(\af)\cap K$, and $\Gamma_{\ell,R}\coloneqq l_R(\Gamma_{R})$ with $l_R$ as in~\eqref{lin}. Set $\Gamma^{red}_{\ell,R}$ to be the restriction of $\Gamma_{\ell,R}$ to $G^{red}_{\ell,R}$. Denote by $\Zhe_{\ell, R}$ (resp. $\Zhe^{red}_{\ell,R}$) the set of all $\Gamma_{\ell,R}$ (resp. $\Gamma^{red}_{\ell,R}$) as we let $K$ vary over the neat compacts with the listed properties.

\subsubsection{The local system on $X(\Gamma^{red}_{\ell,R})$} 

For any congruence subgroup $\Gamma^{red}_{\ell,R} \in\Zhe^{red}_{\ell,R}$ the image of $\Gamma^{red}_{\ell,R}$ under the natural map to
$G^{red}_{\ell,R}(\Qp)$ is contained in $G^{red}_{\ell,R}(\Zp)$.

\begin{lem}\label{factorsp}  The action of $\G^{red}_{\ell,R}$ on $Tor_c^\CO(\boldsymbol{L}^{b}(\cdot, \bar{\pi}_1^*\CJ_{K'}),k)$ factors through the projection on $G^{red}_{\ell,R}(\FF_p)$.
\end{lem}
\begin{proof}  Since the tensor category of automorphic vector bundles  on $M_F$ is generated by the relative cotangent bundle of $A_F$ over $M_F$, 
this follows directly from the description of the latter in \cite[(1.3.2.8)]{Lan3}, given that the action of $\G^{red}_{\ell,R}$ is linear on the factor denoted $X$ in that formula.  Although the reference in \cite{Lan3} is asserted for the compactification in characteristic zero, the statements remain true in mixed characteristic as long as we are in a situation
of good reduction, which is the case throughout this paper.
\end{proof}

\begin{cor}\label{subqu}  Every irreducible constituent of the local system 
$$Tor_c^\CO(\boldsymbol{L}^{b}(\cdot, \bar{\pi}_1^*\CJ_{K'}),k)$$ 
is a subquotient of the local system attached to an algebraic representation of $G^{red}_{\ell,R}$.
\end{cor}
\begin{proof}  This follows from Lemma \ref{factorsp} and
Steinberg's theorem \cite{S63} that every irreducible representation of $G^{red}_{\ell,R}(\FF_p)$ is a subquotient of the restriction to the $\FF_p$ points of an
algebraic representation.
\end{proof}

\subsection{Proof of Claim \ref{Hecke2}}\label{Hecke2proof}   Let $\nu':  T_{T,K} \ra k$ be a character as in Claim \ref{Hecke1}.
It follows from Corollary \ref{subqu} that there is an irreducible algebraic representation $r_\CW:  G_{\ell,R}(\FF_p) \ra \Aut(\CW)$ over $\FF_p$, as in the discussion preceding  Theorem \ref{GaloisrepsR},1 such that $\nu'$ is realized in an $r_\CW$-isotypic subquotient of
$$H^{a}(\mathfrak{N}_{\Sigma}(R),Tor_c^\CO(\boldsymbol{L}^{b}(\cdot, \bar{\pi}_1^*\CJ_{K'}),k))$$
for an automorphic vector bundle $\CJ_{K'}$ on $\bbS_{K_R}(G_{h,R},X(R))$ and for some $a, b, c$.   Part (a) of Claim \ref{Hecke2} then follows
by filtering $\CJ_{K'}$ by irreducible subquotients, one of which will be the $\CF^R$ of the statement of Claim \ref{Hecke2}.  Part (b) is then standard.
More precisely, it comes down to the following observation.  Let $\ell \neq p$ be an unramified prime for $K$ and $\alpha$ be an irreducible spherical representation of $L_R(\QQ_\ell)$ with coefficients in $k$; let $v_R \in \alpha$ be a non-zero spherical vector.  Let $I(\alpha)$ denote the unnormalized induction of the pullback of $\alpha$ to $R(\QQ_\ell)$ to $G(\QQ_\ell)$, and let 
$v \in I(\alpha)$ be a (non-zero) spherical vector.  Let $T(G,K_\ell)$ and $T_\ell(L_R,K_{R,\ell})$ denote the local spherical Hecke algebras of $G(\QQ_\ell)$ and
$L_R(\QQ_\ell)$, respectively, relative to the indicated hyperspecial maximal compact subgroups.
Then $T(G,K_\ell)$ (resp. $T_\ell(L_R,K_{R,\ell})$) acts on $v$ (resp. on $v_R$) through a character
$$\nu_{I(\alpha)}:  T(G,K_\ell) \ra k \text{(resp.  $\nu_\alpha:  T_\ell(L_R,K_{R,\ell}) \ra k$}$$
and
$$\nu_{I(\alpha)} = \nu_\alpha \circ s_R,$$
where $s_R$ is the unnormalized partial Satake transform, as above.  See \cite[Proposition 3.8]{NT} for the analogous case of topological cohomology, and see 
\cite[Lemma 2.3.2]{H84} for the analogous global argument for holomorphic modular forms.

%{\color{red}  Explain how the restriction commutes with partial Satake transform.}

%%%%%%%%%%%%%%%%%%%%%%%%%%%%%%%%%%%%%%%%%%%%%%%%
%%%%%%%%%%%%%%%%%%%%%%%%%%%%%%%%%%%%%%%%%%%%%%%%%
%\section{Boundary cohomology}  

%\begin{thm}  Let $k$ be a finite field and let $\sigma$ be an irreducible representation of $L_P(\af)$, with coefficients in $k$,
%that occurs as a subquotient of the cohomology of $X_U(L_P)$ with coefficients in a local coefficient system $W$.  Then there is an
%$n$-dimensional  representation $\rho(\sigma)$ of $Gal(\Qbar/F)$, with values in $k$, that is associated to $\sigma$ in the sense
%that the characteristic polynomials of $Frob_v$ at unramified places $v$ satisfy \eqref{charpoly}.  Moreover, the representation
%$\rho(\sigma)$ is {\bf reducible}.
%\end{thm}
%\begin{proof}
%\end{proof}

%\vfill
%\pagebreak

\appendix
\section{Adelic boundary components}
\subsection{Parabolic subgroups}\label{parab}

 In what follows, $R$ denotes a standard rational proper parabolic
subgroup of $G$.  We let $r(R)$ denote the parabolic rank of $R$; thus $r(R) = 1$ if and only if
$R$ is a maximal proper parabolic.   Let $L_R$ denote the Levi quotient of $R$.  It is a connected reductive
group that admits a factorization
\begin{equation}\label{LR}
L_R = G_{h,R} \cdot G_{\ell,R}
\end{equation}
where $G_{h,R} = G(V_R)$ is the similitude group of a hermitian vector space $V_R$ over $F$, of dimension $n - 2m(R)$, for
some $1 \leq m(R) \leq \frac{n}{2}$, and 
$$G_{\ell,R} = \prod_{i = 1}^{r(R)} GL(m_i(R))_F, \sum_i m_i(R) = m(R).$$
The factorization is defined by choosing a flag 
$$0 = A_0 \subset A_1 \subset A_2 \subset \dots \subset A_{r(R)}$$
of totally isotropic subspaces of $V$, all assumed to be defined over $\CO$.  Then $GL(m_i(R))$ is identified with the group 
scheme $GL(A_i/A_{i-1})$. over $Spec(\CO)$.  In particular, there is a surjective homomorphism
\begin{equation}\label{lin}  \ell_R:  L_R \ra G_{\ell,R} 
\end{equation}
whose kernel is isomorphic to the unitary similitude group $G(V_R)$, where $V_R$ is the quotient of $V/A_{r(R)}$ by the null space of the induced
hermitian form.  
%\end{rmk}

If $P$ is a  standard rational maximal parabolic subgroup then $r(P) = 1$, $G_{\ell,P} = R_{F/\QQ}GL(m(P))_F$, and $G_{h,P}$ is the similitude group of
a hermitian vector space of dimension $n - 2m(P)$.  Intersection with $G_{\ell,P}$ defines a bijection
between the set of standard rational parabolics $R$ with $G_{h,R} = G_{h,P}$ and the set of standard rational parabolics $R_{\ell,P} \subset G_{\ell,P}$.
In this way we obtain
a canonical map from the set of all standard rational proper parabolic subgroups to the set of (standard rational) maximal parabolic subgroups:  any $R$ is
contained in a unique maximal $P$ with $G_{h,R} = G_{h,P}$, and we let $P(R)$ denote that $P$; we say that $R$ is {\it subordinate} to $P$.

The  (standard rational) maximal parabolic subgroups are totally ordered:  $P < P'$ if and only if $m(P) > m(P')$, which is true if and only if
$\dim G_{h,P} < \dim G_{h,P'}$.   Moreover, $P < P'$ if and only if the standard boundary component $F$ of $X$ stabilized by $P$, as we discuss in \S \ref{termin} below, is smaller than (in fact, is a boundary component of) the boundary component $F'$ stabilized by $P'$.

\subsection{Parametrization of adelic boundary components}\label{param_boundary}

The subgroup $G_{h,P}$ is defined differently than in \cite{HZ94,HZ94b,HZ01}.  In those references the split component of $G_{\ell,P}$ was included as a
subgroup of $G_{h,P}$ in order to account for the mixed Hodge structure on the boundary cohomology; an extra factor of $\Gm$ is needed to define the
appropriate Shimura datum.  The twist does not change the algebraic structure of the Shimura variety attached to $G_{h,P}$ but it does affect the weights
of the variations of Hodge structure on the local systems, and this is reflected in the Tate twists that appear in \eqref{Tatetwist}.  

The factorization \eqref{LR} is a {\bf direct product}.  This implies that the groups denoted $\Delta_{1,R}$ and $\Delta_{0,R}$ that
appear in the formula \cite[(3.2.8)]{HZ01} are trivial.  It also implies that the combinatorial structure of the toroidal boundary bears a simple relation
to the topology of the Borel-Serre compactifications of the locally symmetric spaces attached to $G_{\ell,P}$, as $P$ varies.

We begin by recalling the Borel-Serre compactifications for 
$$G_{\ell,P} = R_{F/\QQ}GL(m(P))_F.$$  Good references for this are \cite[\S 4]{HR20} (though this only treats totally real fields the method is general), \cite[\S 1.3.9]{R}, and \cite{NT}, who consider torsion coefficients, as well as the unpublished book \cite{Har}.   Fix once and for all a rational minimal parabolic
$P_0 \subset G$ with Levi decomposition $P_0 = L_0 \cdot U_0$.  We consider open compact subgroups
$K \subset G(\af)$ with the following properties:
\begin{itemize}
\item $K = K^p \times K_p$ where $K_p$ is a fixed hyperspecial maximal compact subgroup of $G(\Qp)$.
\item For any standard rational parabolic $R \supset P_0$, with $R = L_R \cdot U_R$, $L_R \supset L_0$, we have 
$$K^p \cap R(\afp) = (K^p \cap L_R(\afp))\cdot (K^p \cap U_R(\afp);$$
Let 
$$K_R = K \cap R(\af) = K_p\cap R(\Qp) \times K^p \cap R(\afp).$$
\item $K(L_R) := K\cap L_R(\af) = K_{h}(R) \times K_{\ell}(R)$ where 
$$K_{h}(R) = K(L_R)\cap G_{h,R}(\af);  K_{\ell}(R) =  K(R)\cap G_{\ell,R}(\af).$$
\item $K$ is neat.
\end{itemize}
Such $K^p$ are cofinal in the set of all open compact subgroups of $G(\afp)$.  

We define the locally symmetric space
$$X_K(G_{\ell,P}) = G_{\ell,P}(\QQ)\backslash G_{\ell,P}(\af)/K_{\infty,\ell,P}\cdot K_\ell(P)$$
for a subgroup $K_{\infty,\ell,P}$ containing the center $Z_{G_{\ell,P}}(\RR)$ of $G_{\ell,P}(\RR)$ and 
maximal compact modulo $Z_{G_{\ell,P}}(\RR)$.  We choose maximal compact (mod center) subgroups $K_{\infty,\ell,R}$ 
compatibly with a fixed maximal compact (mod center) subgroup $K_\infty \subset G(\RR)$.

The inclusion of $X_K(G_{\ell,P})$ in its Borel-Serre compactification 
$$X_K(G_{\ell,P}) \subset X_K(G_{\ell,P})^{\rm BS}$$
 is a homotopy equivalence.  The complement
 $$\partial_{\ell,P,K} = X_K(G_{\ell,P})^{\rm BS} ~\setminus ~ X_K(G_{\ell,P})$$
 is the disjoint union of locally closed strata indexed by the rational standard parabolics $R$ of $G$ subordinate to $P$:  
 $$\partial_{\ell,P,K} = \coprod_{P(R) = P} \partial_{R,K}.$$
 For each such $R$, or equivalently for each standard rational parabolic $R_{\ell,P} \subset G_{\ell,P}$, we have
 $$\partial_{R,K} = R_{\ell,P}(\QQ)\cdot U_{R_{\ell,P}}(\af)\backslash R_{\ell,P}(\RR) \times G_{\ell,P}(\af)/K_{\infty,\ell,R}\cdot K_\ell(P),$$
 where $R_{\ell,P}(\QQ)$ acts diagonally on the product, $U_{R_{\ell,P}}(\af)$ and $K_\ell(P)$ act on the factor $G_{\ell,P}(\af)$, and 
 $K_{\infty,\ell,R}$ acts on the right on $R_{\ell,P}(\RR)$.

\section*{Acknowledgements}  We thank Wushi Goldring, Vincent Pilloni, Jack Thorne, Najmuddin Fakhruddin, and Chandrashekhar Khare for answering our questions
about aspects of their work that we consulted in the course of writing this paper, and for helping to bring us up to date on the relevant literature.  In particular, we thank Pilloni for pointing to two proofs in his papers of commutativity of Hecke algebras acting on integral coherent cohomology.
The authors are especially grateful to Kai-Wen Lan, who patiently explained where we could find specific references in his papers for
the many results without which nothing in this paper would have been possible.  Finally, we thank the anonymous referee for a careful reading of an earlier version of this paper.

\vfill
\pagebreak

\end{document}